\documentclass[10pt]{article}

\usepackage{amscd,amsmath, amssymb, fancyhdr, graphicx, url}

\numberwithin{equation}{section}

\newcommand{\version}{version 1.1,\ \ Jan 02, 2014}

\makeatletter
\def\x@arrow{\DOTSB\Relbar}
\def\xlongrightarrowfill@{\arrowfill@\relbar\relbar\longrightarrow}
\newcommand{\xlongrightarrow}[2][]{%
        \ext@arrow 0099\xlongrightarrowfill@{#1}{#2}}
\makeatother

\def\eqref#1{(\ref{#1})}
\newcommand{\goth}{\mathfrak}

\newcommand{\arrow}{{\:\longrightarrow\:}}
\newcommand{\Z}{{\Bbb Z}}
\def\C{{\Bbb C}}
\def\P{{\Bbb P}}

\newcommand{\R}{{\Bbb R}}
\newcommand{\Q}{{\Bbb Q}}

\def\1{\sqrt{-1}\:}
\newcommand{\restrict}[1]{{\left|_{{\phantom{|}\!\!}_{#1}}\right.}}
\newcommand{\cntrct}                
{\hspace{2pt}\raisebox{1pt}{\text{$\lrcorner$}}\hspace{2pt}}

\renewcommand{\tilde}{\widetilde}
\renewcommand{\bar}{\overline}
\renewcommand{\phi}{\varphi}
\renewcommand{\epsilon}{\varepsilon}
\renewcommand{\geq}{\geqslant}
\renewcommand{\leq}{\leqslant}

\newcommand{\Teich}{\operatorname{\sf Teich}}
\newcommand{\Univ}{\operatorname{\sf Univ}}
\newcommand{\Comp}{\operatorname{\sf Comp}}

\newcommand{\Per}{\operatorname{\sf Per}}
\newcommand{\Perspace}{\operatorname{{\Bbb P}\sf er}}

\newcommand{\im}{\operatorname{im}}

\newcommand{\Kah}{\operatorname{Kah}}

\newcommand{\Bir}{\operatorname{Bir}}

\newcommand{\Id}{\operatorname{Id}}

\newcommand{\Hom}{\operatorname{Hom}}

\newcommand{\Pic}{\operatorname{Pic}}
\newcommand{\Pos}{\operatorname{Pos}}

\newcommand{\Aut}{\operatorname{Aut}}
\newcommand{\Mor}{\operatorname{Mor}}
\newcommand{\Mon}{\operatorname{\sf Mon}}

\newcommand{\Diff}{\operatorname{\sf Diff}}

\newcommand{\codim}{\operatorname{codim}}

\newcommand{\Hdg}{\operatorname{\sf Hdg}}

\newcommand{\rk}{\operatorname{rk}}
\newcommand{\Def}{\operatorname{Def}}
\newcommand{\Tw}{\operatorname{Tw}}

\renewcommand{\Re}{\operatorname{Re}}
\renewcommand{\Im}{\operatorname{Im}}

\newcounter{Mycounter}[section]
\newcounter{lemma}[section]
\newcounter{claim}[section]
\renewcommand{\theclaim}{{Claim \thesection.\arabic{claim}}}
\newcommand{\claim}{%
    \setcounter{claim}{\value{Mycounter}}
    \refstepcounter{claim}
    \stepcounter{Mycounter}
    {\noindent \bf \theclaim:\ }}

\newcounter{sublemma}[section]
\newcounter{corollary}[section]
\renewcommand{\thecorollary}{{Corollary \thesection.\arabic{corollary}}}
\newcommand{\corollary}{%
    \setcounter{corollary}{\value{Mycounter}}
    \refstepcounter{corollary}
    \stepcounter{Mycounter}
    {\noindent \bf \thecorollary:\ }}

\newcounter{theorem}[section]
\renewcommand{\thetheorem}{{Theorem \thesection.\arabic{theorem}}}
\newcommand{\theorem}{%
    \setcounter{theorem}{\value{Mycounter}}
    \refstepcounter{theorem}
    \stepcounter{Mycounter}
    {\noindent \bf \thetheorem:\ }}

\newcounter{conjecture}[section]
\renewcommand{\theconjecture}{{Conjecture \thesection.\arabic{conjecture}}}
\newcommand{\conjecture}{%
    \setcounter{conjecture}{\value{Mycounter}}
    \refstepcounter{conjecture}
    \stepcounter{Mycounter}
    {\noindent \bf \theconjecture:\ }}

\newcounter{proposition}[section]
\renewcommand{\theproposition}
      {{Proposition \thesection.\arabic{proposition}}}
\newcommand{\proposition}{%
    \setcounter{proposition}{\value{Mycounter}}
    \refstepcounter{proposition}
    \stepcounter{Mycounter}
    {\noindent \bf \theproposition:\ }}

\newcounter{definition}[section]
\renewcommand{\thedefinition}
      {{Definition~\thesection.\arabic{definition}}}
\newcommand{\definition}{%
    \setcounter{definition}{\value{Mycounter}}
    \refstepcounter{definition}
    \stepcounter{Mycounter}
    {\noindent \bf \thedefinition:\ }}

\newcounter{example}[section]
\renewcommand{\theexample}{{Example \thesection.\arabic{example}}}
\newcommand{\example}{%
    \setcounter{example}{\value{Mycounter}}
    \refstepcounter{example}
    \stepcounter{Mycounter}
    {\noindent \bf \theexample:\ }}

\newcounter{remark}[section]
\renewcommand{\theremark}{{Remark \thesection.\arabic{remark}}}
\newcommand{\remark}{%
    \setcounter{remark}{\value{Mycounter}}
    \refstepcounter{remark}
    \stepcounter{Mycounter}
    {\noindent \bf \theremark:\ }}

\newcounter{problem}[section]
\newcounter{question}[section]
\makeatletter

\@addtoreset{equation}{section} \@addtoreset{footnote}{section}
\makeatother

\def\blacksquare{\hbox{\vrule width 5pt height 5pt depth 0pt}}
\def\endproof{\blacksquare}

\begin{document}

\begin{center}
{\LARGE\bf
Rational curves on hyperk\"ahler manifolds\\[4mm]
}

Ekaterina Amerik\footnote{Partially supported by AG Laboratory NRU-HSE, 
RF government grant, ag. 11.G34.31.0023, and the NRU HSE
Academic Fund Program for 2013-2014, research grant No. 12-01-0107.},  
Misha
Verbitsky\footnote{Partially supported by RFBR grants
 12-01-00944-Á,  NRU-HSE 
Academic Fund Program in 2013-2014, research grant
12-01-0179, and
AG Laboratory NRU-HSE, RF government grant, ag. 11.G34.31.0023.


}

\end{center}

{\small \hspace{0.15\linewidth}
\begin{minipage}[t]{0.7\linewidth}
{\bf Abstract} \\
Let $M$ be an irreducible holomorphically symplectic manifold. 
We show that all faces of the K\"ahler cone of $M$
are hyperplanes $H_i$ orthogonal to certain homology classes,
called monodromy birationally minimal (MBM) 
classes. Moreover, the K\"ahler cone is a
connected component of a complement of the positive 
cone to the union of all $H_i$. 
We provide several characterizations of the MBM-classes. We show the
invariance of MBM property by deformations, as long as the class in question
stays of type $(1,1)$. For hyperk\"ahler 
manifolds with Picard group generated by a negative class 
$z$, we prove that $\pm z$ is $\Q$-effective
if and only if it is an MBM class. We also prove some results towards
the Morrison-Kawamata cone conjecture for hyperk\"ahler manifolds. 
\end{minipage}
}

\tableofcontents


\section{Introduction}


\subsection{K\"ahler cone and MBM classes: an introduction}

Let $M$ be a hyperk\"ahler manifold, that is, a compact,
holomorphically symplectic K\"ahler manifold. We assume that
$\pi_1(M)=0$ and $H^{2,0}(M)=\C$ (the general case reduces to this by 
\ref{_Bogo_deco_Theorem_}).
In this paper we give a description of the K\"ahler cone
of $M$ in terms of a set of cohomology classes $S\subset H^2(M,\Z)$ 
called {\bf MBM classes} (\ref{mbm}). This set is of topological
nature, that is, it depends only on the deformation type of $M$.

It is known since \cite{_Mori:threefolds_} that the K\"ahler cone
of a Fano manifold is polyhedral. Each of its finitely many 
faces is formed by 
classes vanishing on a certain rational curve: one says that the numerical
class of such a curve
generates an extremal ray of the Mori cone. The notion of an extremal ray
also makes sense for projective manifolds which are not Fano: the number of
extremal rays then does not have to be finite. However, they are
discrete in the half-space where the canonical class restricts negatively.

Huybrechts \cite{_Huybrechts:cone_} and Boucksom \cite{_Boucksom-cone_}
have studied the K\"ahler cone of hyperk\"ahler manifolds (not necessarily
algebraic). They have proved that the K\"ahler classes are exactly those
positive classes (i.e. classes with positive Beauville-Bogomolov-Fujiki 
square; see \ref{_posi_cone_Definition_}) which restrict 
positively to all rational curves (\ref{Boucksom-cone}).

Our work puts these results in a deformation-invariant setting. 
Let $\Pos\subset H^{1,1}(M)$ be the positive cone, and $S(I)$ the set of
all MBM classes which are of type (1,1) on 
$M$ with its given complex structure $I$. Then the K\"ahler cone
is a connected component of $\Pos\backslash S(I)^\bot$, 
where $S(I)^\bot$ is the union of all orthogonal complements
to all $z\in S(I)$.

We describe the MBM classes in terms of the minimal curves
on deformations of $M$, birational maps and the monodromy group
action (Sections \ref{_Defo_rac_Section_}, 5),
and formulate a finiteness conjecture (\ref{_minim_curve_Conjecture_})
claiming that primitive integral MBM classes have bounded square. We deduce
the Morrison-Kawamata cone conjecture from this
conjecture. For deformations
of the Hilbert scheme of points on a K3 surface 
our finiteness conjecture is known
(\cite[Proposition 2]{bayer-hass-tschi}).
This gives a proof of Morrison-Kawamata cone conjecture
for deformations of the Hilbert scheme of K3
(\ref{_M_K_bounded_Theorem_}).
Its proof is independently obtained by Markman and Yoshioka using
different methods (forthcoming).

\subsection{Hyperk\"ahler manifolds}

\definition
A {\bf hyperk\"ahler manifold}
is a compact, K\"ahler, holomorphically symplectic manifold.

\hfill

\definition
A hyperk\"ahler manifold $M$ is called
{\bf simple} if $\pi_1(M)=0$, $H^{2,0}(M)=\C$.

\hfill

This definition is motivated by the following theorem
of Bogomolov. 

\hfill

\theorem \label{_Bogo_deco_Theorem_}
(\cite{_Bogomolov:decompo_})
Any hyperk\"ahler manifold admits a finite covering
which is a product of a torus and several 
simple hyperk\"ahler manifolds.
\endproof

\hfill

\remark
Further on, we shall assume (sometimes, implicitly) that
all hyperk\"ahler manifolds we consider are simple.

\hfill

\remark A hyperk\"ahler manifold naturally possesses a whole 2-sphere of 
complex structures (see Section 2). We shall use the notation $(M,I)$ or
$M_I$ to stress that a particular complex structure is chosen, and $M$
to discuss the topological properties (or simply when there is no risque
of confusion).

\hfill

The Bogomolov-Beauville-Fujiki form was
defined in \cite{_Bogomolov:defo_} and 
\cite{_Beauville_},
but it is easiest to describe it using the
Fujiki theorem, proved in \cite{_Fujiki:HK_}.

\hfill

\theorem\label{_Fujiki_Theorem_}
(Fujiki)
Let $M$ be a simple hyperk\"ahler manifold,
$\eta\in H^2(M)$, and $n=\frac 1 2 \dim M$. 
Then $\int_M \eta^{2n}=c q(\eta,\eta)^n$,
where $q$ is a primitive integer quadratic form on $H^2(M,\Z)$,
and $c>0$ an integer. \endproof

\hfill

\remark 
Fujiki formula (\ref{_Fujiki_Theorem_}) 
determines the form $q$ uniquely up to a sign.
For odd $n$, the sign is unambiguously determined as well.
For even $n$, one needs the following explicit
formula, which is due to Bogomolov and Beauville.
\begin{equation}\label{_BBF_expli_Equation_}
\begin{aligned}  \lambda q(\eta,\eta) &=
   \int_X \eta\wedge\eta  \wedge \Omega^{n-1}
   \wedge \bar \Omega^{n-1} -\\
 &-\frac {n-1}{n}\left(\int_X \eta \wedge \Omega^{n-1}\wedge \bar
   \Omega^{n}\right) \left(\int_X \eta \wedge \Omega^{n}\wedge \bar \Omega^{n-1}\right)
\end{aligned}
\end{equation}
where $\Omega$ is the holomorphic symplectic form, and 
$\lambda>0$.

\hfill

\definition\label{_posi_cone_Definition_}
A cohomology class $\eta \in H^{1,1}_\R(M,I)$ is called
{\bf negative} if $q(\eta,\eta)<0$, and {\bf positive}
if $q(\eta,\eta)>0$. Since the signature of $q$ on $H^{1,1}(M,I)$
is $(1, b_2-3)$, the set of positive vectors is disconnected.
{\bf The positive cone} $\Pos(M,I)$ is the connected component of the set
$\{\eta\in H^{1,1}_\R(M,I)\ \ |\ \ q(\eta,\eta)>0\}$ which contains
the classes of the K\"ahler forms. It is easy to check that the positive
cone is convex.

\hfill

The following theorem is crucial for our work.

\hfill

\theorem\label{Boucksom-cone} (Huybrechts, Boucksom; see 
\cite{_Boucksom-cone_})
Let $(M,I)$ be a simple hyperk\"ahler manifold. The K\"ahler cone $\Kah(M,I)$
can be described as follows:
$$\Kah(M,I)=\{\alpha\in \Pos(M,I)|\alpha\cdot C>0\ \forall C\in RC(M,I)\}$$
where $RC(M,I)$ denotes the set of classes of rational curves on $(M,I)$.
\endproof

\subsection{Main results}

\remark
Let $(M,I)$ be a hyperk\"ahler manifold,
and $\phi:\; (M,I)\dashrightarrow (M,I')$ a bimeromorphic
(also called ``birational'') map to another hyperk\"ahler
manifold (note that by a result of Huybrechts 
\cite{_Huybrechts:basic_}, birational hyperk\"ahler
manifolds are deformation equivalent; that is why there is the same
letter $M$ used for the source and the target of $\phi$).
 Since the canonical bundle of $(M,I)$ and $(M,I')$
is trivial, $\phi$ is an isomorphism in codimension 1 (see for example
\cite{_Huybrechts:basic_}, Lemma 2.6). 
This allows one to identify $H^2(M,I)$ and $H^2(M,I')$.
Clearly, this identification is compatible with the
Hodge structure. Further on, we call $(M,I')$
``a birational model'' for $(M,I)$, and identify
$H^2(M)$ for all birational models.

\hfill


\definition
Let $M$ be a hyperk\"ahler manifold. The 
{\bf monodromy group} of $M$ is a subgroup of $GL(H^2(M,\Z))$
generated by monodromy transforms for all Gauss-Manin local systems.
This group can also be characterized in terms of the
mapping class group action (\ref{_monodro_defi_Remark_}).

\hfill

\definition\label{extremal}
Let $(M,I)$ be a hyperk\"ahler manifold.
A rational homology class
$z\in H_{1,1}(M,I)$ is called {\bf $\Q$-effective}
if $N z$ can be represented as a homology class of a curve,
for some $N \in \Z^{>0}$, and {\bf extremal} 
if for any $\Q$-effective homology classes $z_1, z_2\in H_{1,1}(M,I)$
satisfying $z_1+z_2=z$, the classes $z_1, z_2$ are proportional.

\hfill

In the projective case, a negative extremal class is $\Q$-effective and 
some multiple of it is represented by a rational curve. Moreover the
negative part of the cone of $\Q$-effective classes is locally rational 
polyhedral. This is shown by a version of Mori theory adapted to the case
of hyperk\"ahler manifolds (see \cite{HT-GAFA}).  


The Beauville-Bogomolov-Fujiki form allows one to identify $H^2(M,\Q)$ and
$H_2(M,\Q)$. More precisely, it provides an embedding $H^2(M,\Z)\to
H_2(M,\Z)$ which is not an isomorphism (indeed $q$ is not necessarily
unimodular) but becomes an isomorphism after tensoring with $\Q$.
We thus can talk of extremal classes in $H^{1,1}(M, \Q)$, meaning
that the corresponding classes in $H_{1,1}(M, \Q)$ are extremal.
In general, we shall switch from the second homology to the second 
cohomology as it suits us and transfer all definitions from one
situation to the other one by means of the BBF form, with one
important exception, namely, the notion of  effectiveness. Indeed,
an effective class in $H_{1,1}(M, \Q)$ is a class of a curve, whereas
an effective class in $H^{1,1}(M, \Q)$ is a class of a hypersurface.

We also shall extend the BBF form to $H_2(M,\Z)$ as a rational-valued
quad\-ra\-tic form, and to $H_2(M,\Q)$, without further notice. 

The following property, with which we shall work in this paper, looks 
stronger
than extremality modulo monodromy and birational equivalence, but is
equivalent to it whenever the negative part of the Mori cone is locally
rational polyhedral.

Recall that {\bf a face} of a convex cone in a vector space $V$
is an intersection of its boundary and a hyperplane which 
has non-empty interior (\ref{_faces_Definition_}).

\hfill

\definition\label{mbm}
A non-zero negative rational homology class
$z\in H^{1,1}(M,I)$ is called {\bf monodromy birationally minimal} (MBM)
if for some isometry $\gamma\in O(H^2(M,\Z))$ 
belonging to the monodromy group,
 $\gamma(z)^{\bot}\subset H^{1,1}(M,I)$ contains a face 
of the K\"ahler cone of one of birational
models $(M,I')$ of $(M,I)$.

\hfill

The definition has an obvious counterpart for homology classes of type 
$(1,1)$, where one replaces the orthogonality with respect to the BBF form
by the usual duality between homology and cohomology, given 
by integration of forms over cycles.

\hfill

\remark A face of $\Kah(M,I')$ is, by definition, of maximal dimension 
$h^{1,1}(M,I')-1$. So the definition of $z$ being MBM means that 
$\gamma(z)^{\bot}\cap\partial\Kah(M,I')$ contains an open subset of 
$\gamma(z)^{\bot}$.

\hfill

The point of \ref{mbm} is, roughly, as follows. As we 
shall see, rational curves have nice local deformation properties when they
are minimal (extremal). However, globally a deformation of an extremal
rational curve on a variety $(M,I)$ does not have to remain extremal 
on a deformation $(M,I')$. One can only hope to show the 
deformation-invariance
of the property of $z$ being extremal (as long, of course, as it remains 
of type $(1,1)$) {\bf modulo monodromy and birational
equivalence}. In what follows we shall solve this problem for MBM-classes
and apply it to the study of the K\"ahler cone.

\hfill

The deformation equivalence is especially useful because when $\Pic(M)$ has 
rank one, the somewhat obscure notion of 
monodromy birationally minimal classes becomes much more streamlined.

\hfill

\theorem
Let $(M,I)$ be a hyperk\"ahler manifold, $\rk \Pic(M,I)=1$, and
$z\in H_{1,1}(M,I)$ a non-zero negative rational class. Then $z$ is
monodromy birationally minimal if and only if 
$\pm z$ is $\Q$-effective.

\hfill

{\bf Proof:} See \ref{_MBM_posi_theorem_}. \endproof

\hfill

\definition
Let $(M,I)$ be a hyperk\"ahler manifold.
A negative rational class $z\in H^{1,1}_\Q(M,I)$ is called
{\bf divisorial} if $z=\lambda [D]$ for some effective divisor $D$ and $\lambda\in \Q$.

\hfill

Our first main results concern the deformational invariance of 
these notions.

\hfill

\theorem
Let $(M,I)$ be a hyperk\"ahler manifold, 
$z\in H_{1,1}(M,I)$ an integer homology class, $q(z,z)<0$,
and $I'$ a deformation of $I$ such that $z$ is of type (1,1)
with respect to $I'$. Assume that $z$ is 
monodromy birationally minimal on $(M,I)$. Then
$z$ is monodromy birationally minimal on $(M,I')$.
The property of $z\in H^{1,1}(M,I)$ being divisorial is likewise 
deformation-invariant,
provided that one restricts oneself to the complex structures with
Picard number one (i.e. such that the Picard group is generated by $z$
over $\Q$).

\hfill

{\bf Proof:} See \ref{divisor-invariant}, 
\ref{divisor-pic1} and \ref{MBM-invar-thm}. 
\endproof

\hfill

\remark We expect that the property of $z$ being divisorial is 
deformation-invariant as long as $z$ stays of type $(1,1)$, and plan to
return to this question in a forthcoming paper.

\hfill

The MBM classes can be used to determine the K\"ahler cone of
$(M,I)$ explicitly.

\hfill

\theorem
Let $(M,I)$ be a hyperk\"ahler manifold, 
and $S\subset H_{1,1}(M,I)$ the set of all MBM classes.
Consider the corresponding set of hyperplanes
$S^\bot:=\{W=z^\bot\ \ |\ \ z\in S\}$ in $H^{1,1}(M,I)$.
Then the K\"ahler cone of $(M,I)$ is 
a connected component of $\Pos(M,I)\backslash \cup S^\bot$,
where $\Pos(M,I)$ is the positive cone of $(M,I)$.
Moreover, for any connected component $K$ of 
$\Pos(M,I)\backslash \cup S^\bot$,
there exists $\gamma\in O(H^2(M,\Z))$ in the monodromy group of $M$, and
a hyperk\"ahler manifold $(M,I')$ birationally equivalent 
to $(M,I)$, such that $\gamma(K)$ is the K\"ahler cone of $(M,I')$.

\hfill

{\bf Proof:} See \ref{_KW_chambers_MBM_Theorem_}. \endproof

\hfill

\remark In particular, $z^\bot \cap \Pos(M,I)$ either has dense intersection
with the interior of the K\"ahler chambers (if $z$ is not MBM), or is a union of walls
of those (if $z$ is MBM); that is, there are no 
``barycentric partitions'' in 
the decomposition 
of the positive cone into the K\"ahler chambers.

\centerline{\begin{tabular}{cc}
\includegraphics[width=0.30\textwidth]{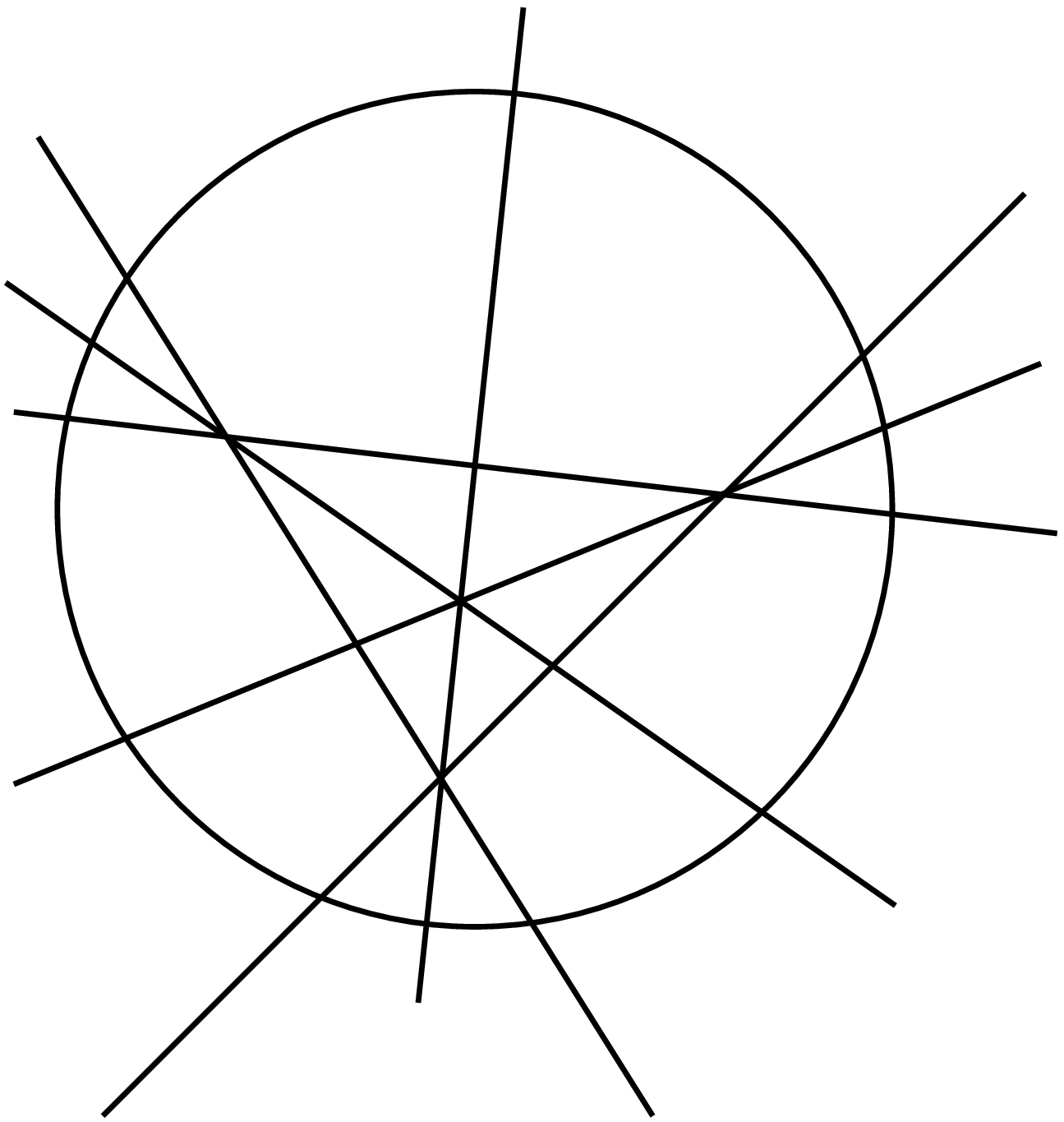}\  &\  \includegraphics[width=0.30\textwidth]{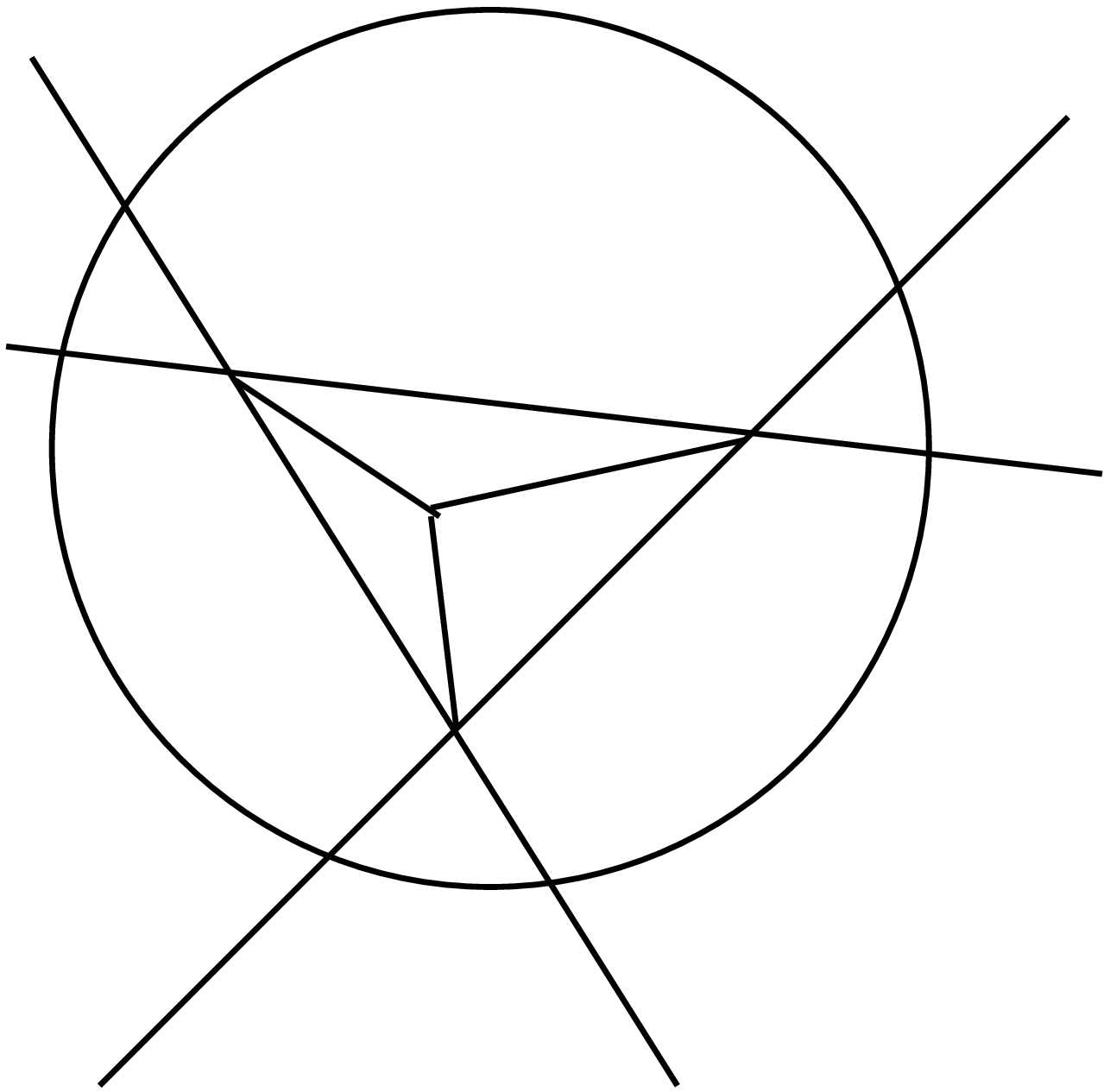} \\
{\bf \scriptsize Allowed partition} & {\bf \scriptsize Prohibited partition} 
\end{tabular}}

\hfill

Finally, we apply this to the Morrison-Kawamata cone conjecture:

\hfill

\theorem\label{_M_K_bounded_Theorem_}
Let $M$ be a simple hyperk\"ahler manifold, and $q$
the Bogomolov-Beauville-Fujiki form. Suppose that 
there exists $C>0$ such that $|q(\eta, \eta)|<C$ for all primitive 
integral MBM classes (or, alternatively, for all extremal rational curves
on all deformations of $M$). Then 
the Morrison-Kawamata cone conjecture holds for $M$: the group $\Aut(M)$
acts on the set of faces of the K\"ahler cone with finitely many orbits.

\hfill

{\bf Proof:} See \ref{morkaw}. \endproof

\hfill

The condition of the theorem is satisfied for manifolds which are deformation
equivalent to the Hilbert scheme of points on a K3 surface, see \cite{bayer-hass-tschi}; 
for such manifolds one therefore obtains a proof of the Morrison-Kawamata 
cone conjecture.


\section{Global Torelli theorem, hyperk\"ahler structures and monodromy group}


In this Section, we recall a number of 
results  about hyperk\"ahler manifolds,
used further on in this paper. For more
details and references, please see 
\cite{_Besse:Einst_Manifo_} and \cite{_V:Torelli_}.

\subsection{Hyperk\"ahler structures and twistor spaces}
\label{_hk_twi_Subsection_}

\definition
Let $(M,g)$ be a Riemannian manifold, and $I,J,K$
endomorphisms of the tangent bundle $TM$ satisfying the
quaternionic relations
\[
I^2=J^2=K^2=IJK=-\Id_{TM}.
\]
The triple $(I,J,K)$ together with
the metric $g$ is called {\bf a hyperk\"ahler structure}
if $I, J$ and $K$ are integrable and K\"ahler with respect to $g$.

Consider the K\"ahler forms $\omega_I, \omega_J, \omega_K$
on $M$:
\[
\omega_I(\cdot, \cdot):= g(\cdot, I\cdot), \ \
\omega_J(\cdot, \cdot):= g(\cdot, J\cdot), \ \
\omega_K(\cdot, \cdot):= g(\cdot, K\cdot).
\]
An elementary linear-algebraic calculation implies
that the 2-form $\Omega:=\omega_J+\1\omega_K$ is of Hodge type $(2,0)$
on $(M,I)$. This form is clearly closed and
non-degenerate, hence it is a holomorphic symplectic form.

In algebraic geometry, the word ``hyperk\"ahler''
is essentially synonymous with ``holomorphically
symplectic'', due to the following theorem, which is
implied by Yau's solution of Calabi conjecture
(\cite{_Besse:Einst_Manifo_}, \cite{_Beauville_}).

\hfill

\theorem\label{_Calabi-Yau_Theorem_}
Let $M$ be a compact, K\"ahler, holomorphically
symplectic manifold, $\omega$ its K\"ahler form, $\dim_\C M =2n$.
Denote by $\Omega$ the holomorphic symplectic form on $M$.
Suppose that $\int_M \omega^{2n}=\int_M (\Re\Omega)^{2n}$.
Then there exists a unique hyperk\"ahler metric $g$ with the same
K\"ahler class as $\omega$, and a unique hyperk\"ahler structure
$(I,J,K,g)$, with $\omega_J = \Re\Omega$, $\omega_K = \Im\Omega$.
\endproof

\hfill

Further on, we shall speak of ``hyperk\"ahler manifolds''
meaning ``holomorphic symplectic manifolds of K\"ahler
type'', and ``hyperk\"ahler structures'' meaning
the quaternionic triples together with a metric. 

\hfill

Every hyperk\"ahler structure induces a whole 2-dimensional
sphere of complex structures on $M$, as follows. 
Consider a triple $a, b, c\in R$, $a^2 + b^2+ c^2=1$,
and let $L:= aI + bJ +cK$ be the corresponging quaternion. 
Quaternionic relations imply immediately that $L^2=-1$,
hence $L$ is an almost complex structure. 
Since $I, J, K$ are K\"ahler, they are parallel with respect
to the Levi-Civita connection. Therefore, $L$ is also parallel.
Any parallel complex structure is integrable, and K\"ahler.
We call such a complex structure $L= aI + bJ +cK$
{\bf a complex structure induced by a hyperk\"ahler structure}.
There is a 2-dimensional holomorphic family of 
induced complex structures, and the total space
of this family is called {\bf the twistor space}
of a hyperk\"ahler manifold. 

\subsection{Global Torelli theorem and monodromy}

\definition
Let $M$ be a compact complex manifold, and 
$\Diff_0(M)$ a connected component of its diffeomorphism group
({\bf the group of isotopies}). Denote by $\Comp$
the space of complex structures on $M$, equipped with
its structure of a Fr\'echet manifold, and let
$\Teich:=\Comp/\Diff_0(M)$. We call 
it {\bf the Teichm\"uller space.}

\hfill

\remark
In many important cases, such as
for manifolds with trivial canonical class (\cite{_Catanese:moduli_}), 
$\Teich$ is a finite-dimensional
complex space; usually it is non-Hausdorff.

\hfill

\definition\label{univfam} The universal family
of complex manifolds over $\Teich$ is defined as $\Univ/\Diff_0(M)$,
where $\Univ$ is the natural universal family over $\Comp$ with its 
Fr\'echet manifold structure. Locally, it is isomorphic to the universal
family over the Kuranishi 
space.\footnote{We are grateful to Claire Voisin for this observation.}

\hfill

\definition The {\bf mapping class group} is $\Gamma=\Diff(M)/\Diff_0(M)$. It
naturally acts on $\Teich$ (the quotient of $\Teich$ by $\Gamma$ may be
viewed as the ``moduli space'' for $M$, but in general it has very bad 
properties; see below).

\hfill

\remark
Let $M$ be a hyperk\"ahler manifold (as usually, we assume
$M$ to be simple). For any $J\in \Teich$,
$(M,J)$ is also a simple hyperk\"ahler manifold, 
because the Hodge numbers are constant in families.
Therefore, $H^{2,0}(M,J)$ is one-dimensional. 

\hfill

\definition
 Let 
\[ \Per:\; \Teich \arrow {\Bbb P}H^2(M, \C)
\]
map $J$ to the line $H^{2,0}(M,J)\in {\Bbb P}H^2(M, \C)$.
The map $\Per$ is 
called {\bf the period map}.

\hfill

\remark
The period map $\Per$ maps $\Teich$ into an open subset of a 
quadric, defined by
\[
\Perspace:= \{l\in {\Bbb P}H^2(M, \C)\ \ | \ \  q(l,l)=0, q(l, \bar l) >0\}.
\]
It is called {\bf the period domain} of $M$.
Indeed, any holomorphic symplectic form $l$
satisfies the relations $q(l,l)=0, q(l, \bar l) >0$,
as follows from \eqref{_BBF_expli_Equation_}.

\hfill

\proposition\label{_period_Grassmann_Proposition_}
The period domain $\Perspace$
is identified with the quotient
$SO(b_2-3,3)/SO(2) \times SO(b_2-3,1)$, which
is the Grassmannian of positive oriented 2-planes in $H^2(M,\R)$.

\hfill

{\bf Proof:} This statement is well known, but we shall
sketch its proof for the reader's convenience.

{\bf Step 1:} Given $l\in {\Bbb P}H^2(M, \C)$, the space
generated by $\Im l, \Re l$ is 2-dimensional, because 
$q(l,l)=0, q(l, \bar l)>0$ implies that $l \cap H^2(M,\R)=0$.

{\bf Step 2:}  This 2-dimensional plane is 
positive, because 
 $q(\Re l, \Re l) = q(l+ \bar l, l+ \bar l) = 2 q(l, \bar l)>0$.

{\bf Step 3:} Conversely, for any 2-dimensional positive
plane  $V\in H^2(M,\R)$, 
the quadric $\{l\in V \otimes_\R \C\ \ | \ \ q(l,l)=0\}$
consists of two lines; a choice of a line is determined by the
orientation.
\endproof

\hfill

\definition
Let $M$ be a topological space. We say that $x, y \in M$
are {\bf non-separable} (denoted by $x\sim y$)
if for any open sets $V\ni x, U\ni y$, $U \cap V\neq \emptyset$.

\hfill

\theorem\label{nonsep-birat} 
(Huybrechts; \cite{_Huybrechts:basic_}).
If two points $I,I'\in \Teich$ are non-separable, then  
there exists a bimeromorphism $(M,I)\dasharrow (M,I')$.
\endproof

\hfill

\remark\label{birat-nonsep}
The converse is not true since many points in the Teichm\"uller 
space correspond to isomorphic manifolds, and they are not always pairwise
non-separable (consider orbits of the
mapping class group action). But it is true that if 
two hyperk\"ahler manifolds are bimeromorphic, then one can find 
non-separable points in the Teichm\"uller space representing them.

\hfill

\definition
The space $\Teich_b:= \Teich\!/\!\!\sim$ is called {\bf the
birational Teichm\"uller space} of $M$.

\hfill

\remark This terminology is slightly misleading since there are 
non-separable points of the Teichm\"uller space which correspond
to biregular, not just birational, complex structures. Even for
K3 surfaces, the Teichm\"uller space is non-Hausdorff. 

\hfill

\theorem \label{_glo_Torelli_Theorem_}
(Global Torelli theorem; \cite{_V:Torelli_})
The period map 
$\Teich_b\stackrel \Per \arrow \Perspace$ is an isomorphism
on each connected component of $\Teich_b$.
\endproof

\hfill

\remark By a result of Huybrechts (\cite{_Huybrechts:finiteness_}), $\Teich$ has only finitely
many connected components. We shall denote by $\Teich_I$ the component 
containing the parameter point for the complex structure $I$, and by 
$\Gamma_I$
the subgroup of the mapping class group $\Gamma$ fixing this component.
Obviously $\Gamma_I$ is of finite index in $\Gamma$.

\hfill

\definition
Let $M$ be a hyperkaehler manifold,
$\Teich_b$ its birational Teichm\"uller space,
and $\Gamma$ the {\bf mapping class group} $\Diff (M)/\Diff_0(M)$.
The quotient $\Teich_b/\Gamma$ is called
{\bf the birational moduli space} of $M$.
Its points are in bijective correspondence with the
complex structures of hyperk\"ahler type on $M$
up to a bimeromorphic equivalence.

\hfill

\remark
The word ``space'' in this context is misleading.
In fact, the quotient topology on $\Teich_b/\Gamma$ is extremely
non-Hausdorff, e.g. every two open sets would intersect
(\cite{_Verbitsky:ergodic_}).

\hfill

The Global Torelli theorem can be stated
as a result about the birational moduli space.

\hfill

\theorem\label{_moduli_monodro_Theorem_}
(\cite[Theorem 7.2, Remark 7.4, Theorem 3.5]{_V:Torelli_})
Let $(M,I)$ be a hyperk\"ahler manifold, and $W$ 
a connected component of its birational
moduli space. Then $W$ is isomorphic to ${\Perspace}/\Mon_I$,
where ${\Perspace}=SO(b_2-3,3)/SO(2) \times SO(b_2-3,1)$
and $\Mon_I$ is 
an arithmetic group in $O(H^2(M, \R), q)$, called {\bf the
monodromy group} of $(M,I)$. In fact $\Mon_I$ is the image of $\Gamma_I$
in $O(H^2(M, \R), q)$.  
\endproof

\hfill

\remark\label{_monodro_defi_Remark_}
The monodromy group of $(M,I)$ can be also described
as a subgroup of the group $O(H^2(M, \Z), q)$
generated by monodromy transform maps for 
Gauss-Manin local systems obtained from all
deformations of $(M,I)$ over a complex base
(\cite[Definition 7.1]{_V:Torelli_}). This is 
how this group was originally defined by Markman
(\cite{_Markman:constra_}, \cite{_Markman:survey_}).
The fact that it is of finite index in $O(H^2(M, \Z), q)$
is crucial for the Morrison-Kawamata conjecture, see next section. 

\hfill

\remark 
A caution: usually ``the global Torelli theorem''
is understood as a theorem about Hodge structures.
For K3 surfaces, the Hodge structure on $H^2(M,\Z)$
determines the complex structure. 
For $\dim_\C M >2$, it is false.


\section{Morrison-Kawamata cone conjecture for
  hyperk\"ahler manifolds}


In this section, we
introduce the Morrison-Kawamata cone conjecture,
and give a brief survey of what is known, following
\cite{_Markman:survey_} and \cite{_Totaro:MK_cone_}, with some easy
generalizations.

Originally, this conjecture was stated in \cite{_Morrison:Beyond_}
and proven by Kawamata in \cite{_Kawamata:cone_}  
for Calabi-Yau threefolds admitting a holomorphic fibtation 
over a positive-dimensional base.

\subsection{Morrison-Kawamata cone conjecture}

\definition\label{_faces_Definition_}
Let $M$ be a compact, K\"ahler manifold,
$\Kah\subset H^{1,1}(M,\R)$ the K\"ahler cone,
and $\overline\Kah$ its closure in $H^{1,1}(M,\R)$,
called {\bf the nef cone}. A {\bf face} of the 
K\"ahler cone is the intersection
of the boundary of $\overline\Kah$ and a hyperplane $V\subset H^{1,1}(M,\R)$,
which has non-empty interior.

\hfill

\conjecture (Morrison-Kawamata cone conjecture) \\
Let $M$ be a Calabi-Yau manifold. Then 
the group $\Aut(M)$ of biholomorphic automorphisms of $M$
acts on the set of faces of $\Kah$  with finite number of orbits.

\hfill

This conjecture has a birational version, which has many
important implications. For projective hyperk\"ahler manifolds,
the birational version of Morrison-Kawamata
cone conjecture has been proved by E. Markman in 
\cite{_Markman:survey_}.

\hfill

\definition
Let $M,M'$ be compact complex manifolds.
Define a {\bf  pseudo-isomorphism} $M\dasharrow M'$ as a
birational map which is an isomorphism outside of codimension
$\geq 2$ subsets of $M, M'$.

\hfill

\remark For any pseudo-isomorphic manifolds $M,M'$,
the second cohomologies $H^2(M)$ and $H^2(M')$ are naturally
identified.

\hfill

As we have already remarked, any birational map of hyperk\"ahler
varieties is a pseudo-isomorphism; more generally, this is true
for all varieties with nef canonical class.

\hfill

\definition The
{\bf movable cone}, also known as {\bf  birational nef cone}
is the closure of the union of $\Kah(M')$ for all
$M'$ pseudo-isomorphic to $M$. 
The union of $\Kah(M')$ for all
$M'$ pseudo-isomorphic to $M$ is called {\bf birational K\"ahler cone}.

\hfill

\conjecture (Morrison-Kawamata birational cone conjecture)\\
Let $M$ be a Calabi-Yau manifold. Then 
the group $\Bir(M)$ of birational automorphisms of $M$
acts on the set of faces of its movable cone with finite 
number of orbits.

\hfill

\subsection{Morrison-Kawamata cone conjecture for K3 surfaces}

In this subsection, we prove the Morrison-Kawamata 
cone conjecture for K3 surfaces. Originally it was proven
by Sterk (see \cite{_Sterk:finiteness_}). The proof we are
giving has two advantages: it works in the non-algebraic
situation, and to some extent generalizes to arbitrary dimension.
Notice that the pseudo-isomorphisms of smooth surfaces are
isomorphisms, hence for K3 surfaces both flavours of
Morrison-Kawamata cone conjecture are equivalent.

\hfill

\definition
A cohomology class $\eta\in H^2(M,\Z)$ on a K3 surface
is called {\bf a $(-2)$-class} if
$\int_M\eta\wedge\eta=-2$.

\hfill

\remark
Let $M$ be a K3 surface, and $\eta\in H^{1,1}(M,\Z)$ a $(-2)$-class.
 Then either $\eta$ or $-\eta$ is effective. Indeed,
$\chi(\eta) = 2+\frac{\eta^2}{2}=1$ by Riemann-Roch.

\hfill

The following theorem is well-known.

\hfill

\theorem
Let $M$ be a K3 surface, and $S$ the set of all effective $(-2)$-classes.
Then $\Kah(M)$ is the set of all $v\in \Pos(M)$
such that $q(v, s) >0$ for all $s\in S$.

\hfill

{\bf Proof:} By adjunction formula, a
curve $C$ on a K3 surface satisfies $C^2=2g-2$, where $g$
is its genus, hence the cone of effective curves is
generated by positive classes and effective $(-2)$-classes.
By \cite{_Demailly_Paun_}, the K\"ahler cone is the
intersection of the positive cone and the cone dual to the
cone of effective curves, hence it is $v\in \Pos(M)$
such that $q(v, s) >0$ for all effective $s$.
However, for all $s,v\in \Pos(M)$
the inequality $q(v, s) >0$ automatically
follows from the Hodge index formula.
\endproof

\hfill

\definition
{\bf A Weyl chamber}, or {\bf K\"ahler chamber}, on a K3 surface is a connected 
component
of $\Pos(M)\backslash S^\bot$, where $S^{\bot}$ is the union of
all planes $s^\bot$ for all $(-2)$-classes $s\in S$. 
{\bf The reflection group}
of a K3 surface is the group $W$ generated by reflections with
respect to all $s\in S$.

\hfill

\remark
Clearly, a Weyl chamber is a fundamental domain of $W$,
and $W$ acts transitively on the set of all Weyl chambers.
Moreover, the K\"ahler cone of $M$ is one of its Weyl
chambers.

\hfill

To a certain extent, this is a pattern which is repeated in all
dimensions, as we shall see.

\hfill

The following theorem is well-known, too, but we want to sketch a proof 
in some detail since it is also important for what follows.

\hfill 

\theorem\label{_K3_auto_Theorem_}
Let $M$ be a K3 surface, and $\Aut(M)$ 
the group of all automorphisms of $M$. 
Then $\Aut(M)$ is the group of all 
isometries of $H^2(M,\Z)$ preserving the K\"ahler
chamber $\Kah$ and the Hodge decomposition.

\hfill

{\bf Proof:} First, let us show that the natural map
$\Aut(M)\stackrel\phi\to O(H^2(M,\Z))$ is injective.
Our argument is similar to Buchdahl's in \cite{Buchdahl}. 
Clearly, $\Aut(M)$ acts on $H^2(M,\Z)$;
its kernel is formed by automorphisms preserving all
K\"ahler classes, hence acting as isometries on 
all Calabi-Yau metrics on all deformations of $M$.
Deforming $M$ to a Kummer surface, we find that
this isometry must also induce an isometry on the
underlying 2-torus, acting trivially on cohomology.
Therefore $\ker \phi$ acts as identity on all Kummer surfaces,
and hence on all K3 surfaces as well. 

To describe the image of $\phi$,
we use the Torelli theorem. Let $\gamma\in O(H^2(M,\Z))$
preserve the Hodge decomposition on $(M,I)$ and the
K\"ahler cone. 
Global Torelli theorem affirms that the 
period map $\Per:\; \Teich_b \arrow \Perspace$
for K3 surfaces is an isomorphism. The mapping class
group $\Gamma=\Diff(M)/\Diff_0(M)$ acting on $\Teich$ embeds as  
the orientation-preserving subgroup in $O(H^2(M,\Z))$ (\cite{Borcea}).
We deduce that $\gamma$ comes from an element $\tilde\gamma$ of the mapping class group.
Since $\gamma$ preserves the Hodge decomposition, the 
corresponding point in the period space
$\Per(I)\in \Perspace$ is a fixed point of $\tilde\gamma$.
Therefore, $\tilde\gamma$ exchanges the non-separable points in 
$\Per^{-1}(\Per(I))$. However, as shown by Huybrechts,
two manifolds which correspond to distinct non-separable 
points in $\Teich$ have different K\"ahler cones, so that 
the fibers of
the period map parametrize the K\"ahler (Weyl) chambers
(\cite{_Huybrechts:basic_}). Therefore, $\tilde\gamma$
maps $(M,I)\in \Per^{-1}(\Per(I))$ to itself whenever it fixes the
K\"ahler cone. Hence $\tilde\gamma$, considered as a diffeomorphism of $M$, is an automorphism
of the complex manifold $(M,I)$, and $\gamma$ is in the image of $\phi$ . \endproof

\hfill

Now we can prove the Morrison-Kawamata cone conjecture for K3 surfaces.
Our argument is based on the following general result
about lattices.

\hfill

\theorem\label{_Kneser_orbits_Theorem_}
Let $q$ be an integer-valued quadratic form  on
$\Lambda=\Z^n,\ n\geq 2$ (not necessarily unimodular) such that its kernel
is at most one-dimensional, and 
$O(\Lambda)$ the corresponding group of isometries.
Fix $0\neq r\in \Z$, and let $S_r$ be the set $\{v\in
\Lambda\ \ |\ \ q(v,v)=r\}$. Then $O(\Lambda)$
acts on $S_r$ with finite number of orbits.

\hfill

{\bf Proof:} The non-degenerate case is
\cite[Satz 30.2]{_Kneser:Quadratische_}. In the case of one-dimensional
kernel $L=\ker q$, we write $\Lambda=L\oplus \Lambda_0$ and the 
elements of $\Lambda$ as $a_0+kl$ where $a_0\in \Lambda_0$,
$k\in \Z$ and $l$ is a fixed generator of $L$. There are finitely many
representatives of $S_r\cap \Lambda_0$ under the action of $O(\Lambda_0)$, 
say
$a_0^1,\dots, a_0^m$. For each $j$, take a system of representatives
$k_1,\dots, k_{t_j}$ of $\Z$ modulo the ideal $\Hom(\Lambda_0, L)\cdot a_0^j$.
Then representatives 
of the orbits in $S_r$ are elements of the form $a_0^j+k_il$.  \endproof

\hfill

We shall apply this to study the Picard lattice $H^2(M,\Z)\cap H^{1,1}(M)$. 
Notice that it is non-degenerate when $M$ is projective, and can have
at most one-dimensional kernel when $M$ is arbitrary.

\hfill

\theorem \label{_cone_for_K3_Corollary_}
 Morrison-Kawamata cone conjecture holds for any K3 surface.

\hfill

{\bf Proof, step 1: setting a goal.} Let $\Mon^{\Hdg}(M,I)$ be the group of all
oriented isometries of $H^2(M,\Z)$ preserving the Hodge
decomposition (such isometries are known as {\bf Hodge
  isometries}). We have already remarked that when $M$ is a K3 surface, all these are
monodromy operators, hence the notation.  Since $\Mon^{\Hdg}(M,I)$ acts on the Picard lattice $\Lambda$ as a subgroup 
of index at most two in $O(\Lambda)=O(H^{1,1}(M,\Z))$, 
\ref{_Kneser_orbits_Theorem_} implies that $\Mon^{\Hdg}(M,I)$
acts with finitely many
orbits on the set of $(-2)$-vectors in $\Lambda$.

Our goal in the next steps is to prove that the group $\Mon^{\Hdg}(M,I)$ 
acts on the set of faces
of the Weyl chambers with finitely many orbits.

\hfill

{\bf Step 2: describing a face by a flag.} Recall that an orientation of a 
hyperplane in an oriented  real vector
space is the choice of a ``positive'' normal direction. Then
a face $F$ of a Weyl chamber is determined by the following data:
a hyperplane $P_{s-1}$ it sits on, with the orientation
pointing to the interior of the chamber; a hyperplane $P_{s-2}$ in 
$P_{s-1}$ supporting one of the faces $F_1$ of $F$, together with the
orientation pointing to the interior of the face $F$; an oriented hyperplane
$P_{s-3}$ in $P_{s-2}$ supporting some face $F_2$ of $F_1$; and so forth.
Here $s$ is the dimension of the ambient space $H^{1,1}_{\R}(M)$.
In other words, a face is determined by a full flag of
linear subspaces, oriented step-by-step as above. 

To prove that $\Mon^{\Hdg}(M,I)$ acts on the set of faces with finitely
many orbits, it suffices to prove that it acts with finitely many orbits
on flags $P_{s-1}\supset P_{s-2}\supset \dots \supset P_1$ of the above 
form, orientation
forgotten.

\hfill

{\bf Step 3: bounding squares.} Notice that the hyperplane $P_{s-1}$ is $x^\bot$, where 
$x$ is a $(-2)$-class in $H^{1,1}_{\Z}(M,I)$, and $P_{s-2}=x^\bot\cap y^\bot$ 
where $y$ is another $(-2)$-class in $H^{1,1}_{\Z}(M,I)$. We claim
that $P_{s-2}$ as a hyperplane in $P_{s-1}$ is given by $(y')^\bot$,
where $y'\in x^\bot$ is an integral negative class of square at least $-8$.

Indeed, write $y=\frac{\langle x,y\rangle}{\langle x,x\rangle}x+\tilde{y}$
(orthogonal projection to $x^\bot$). Obviously,  $P_{s-2}$ as a hyperplane in
 $P_{s-1}$ is just the orthogonal to $\tilde{y}\in x^\bot$. The first summand on the right has
non-positive square since it is proportional to a $(-2)$-class $x$.
We claim that the second summand has strictly negative square (this 
square has then to be at least $-2$). Indeed, since
the hyperplanes $x^\bot$ and $y^\bot$ define a K\"ahler chamber, the 
intersection of $x^\bot$ and $y^\bot$ is within the positive cone. But then
the orthogonal to $\tilde{y}$ in $x^\bot$ intersects the positive cone
in $x^\bot$, so that, since the signature of the intersection form restricted
to $x^\bot$ is $(+,-,\dots, -)$, $\tilde{y}$ is
of strictly negative square. To make $\tilde{y}$ integral, it suffices
to replace it by $y'=\langle x,x\rangle \tilde{y}=-2\tilde{y}$; one has
$0>\langle y', y'\rangle \geq -8$.

\hfill

{\bf Step 4: monodromy action.} We have just seen that the intersections of 
$x^\bot$ with the other hyperplanes defining the K\"ahler chambers are given,
in $x^\bot$, as orthogonals to integer vectors of strictly negative bounded 
square. By \ref{_Kneser_orbits_Theorem_}, the orthogonal group 
of $H^{1,1}_{\Z}(M,I)\cap x^\bot$ acts with finitely many orbits on such
vectors. We deduce that the orthogonal group of $\Lambda=H^{1,1}_{\Z}(M,I)$ acts
with finitely many orbits on partial flags $P_{s-1}\supset P_{s-2}$ arising
from the faces of K\"ahler chambers. Iterating the argument of Step 3, we
see that $O(\Lambda)$, and thus also $\Mon^{\Hdg}(M,I)$, acts with finitely many 
orbits on full flags and therefore on the set of faces of K\"ahler chambers.
    
\hfill

{\bf Step 5: conclusion.}
For each pair of faces $F, F'$ of a K\"ahler
cone and $w\in \Mon^{\Hdg}(M,I)$ mapping $F$ to $F'$, $w$ maps $\Kah$ to 
itself
or to an adjoint Weyl chamber $K'$. Then $K'=r(K)$, where
$r$ is the orthogonal reflection in $H^2(M,\Z)$ fixing $F'$.
In the first case, $w\in \Aut(M)$.
In the second case, $rw$ maps $F$ to $F'$ and maps $\Kah$ to itself,
hence $rw\in \Aut(M)$.
\endproof

\subsection{Finiteness of polyhedral tesselations}

The argument used to prove 
\ref{_cone_for_K3_Corollary_}
is valid in an abstract setting, which 
is worth describing here. 

\hfill

Let $V_\Z$ is a torsion-free $\Z$-module
of rank $n+1$, equipped with an integer-valued (but
not necessarily unimodular) quadratic form of signature $(1,n)$, $\Gamma\subset O(V_\Z)$ a finite index subgroup,
and $V=V_\Z\otimes_\Z \R$. Let $S_0\subset V_\Z$ be a finite 
set of negative vectors,
$S:=\Gamma\cdot S_0$ its orbit, and $Z:=\bigcup_{s\in S} s^\bot$
the union of all orthogonal complements to $s$ in $V$.

We associate to $V$ a hyperbolic space ${\Bbb H}={\Bbb P}V^+$
obtained as a projectivization of the set of positive vectors.
Let $\{P_i\}$ be the set of connected components of 
${\Bbb P}V^+\backslash {\Bbb P}Z$. This is a polyhedral
tesselation of the hyperbolic space. 

\hfill

\definition\label{_cut_tessela_Definition_}
We call such a tesselation
{\bf a tesselation 
cut out by the set of hyperplanes orthogonal to $S$}.

\hfill

\definition
Let $S_d$ be an intersection of $n-d$
transversal hyperplanes $s^\bot$, with $s\in S$.
A {\bf $d$-dimensional face} of a tesselation is
a connected component of $S_d\backslash S_{d-1}$.

\hfill

\theorem\label{_tessela_fini_orbits_Theorem_}
Let $\{P_i\}$ be a tesselation obtained in
\ref{_cut_tessela_Definition_}, and $F_d$ the
set of all $d$-dimensional faces of the tesselation.
Then $\Gamma$ acts on $F_d$ with finitely many orbits.

\hfill

{\bf Proof:}
We encode a $d$-dimensional face $S_d^0$ (up to a finite choice of
orientations) by a sequence of hyperplanes
$s_1^\bot$, $s_2^\bot$, ... $s_{n-d}^\bot$ such that their intersection 
supports $S_d^0$, plus an oriented flag as in \ref{_cone_for_K3_Corollary_}.
Using induction, we may assume that up to the action of $\Gamma$ there 
are only finitely many $(d+1)$-dimensional faces. 
The same inductive argument as in \ref{_cone_for_K3_Corollary_},
Step 3 is used to show that the projection of
$s_{n-d}$ to $S_{d+1}:=\bigcap_{i=1}^{n-d-1}s_i^\bot$ has bounded
square after multiplying by the denominator. 
Therefore, there are finitely many orbits of the group
$\Gamma_{S_{d+1}}:=\{\gamma\in \Gamma\ \ |\ \ \gamma(S_{d+1})=S_{d+1}\}$ 
on the set $s_{n-d}^\bot\cap S_{d+1}$ for all $s_{n-d}\in S$, and the reasoning is
the same for the oriented flag.
\endproof

\subsection{Morrison-Kawamata cone conjecture: birational version}




\definition The birational K\"ahler cone of a hyperk\"ahler manifold $M$
is the union of pullbacks of the K\"ahler cones under all birational 
maps $M\dasharrow M'$ where $M'$ is also a hyperk\"ahler 
manifold.

\hfill

\remark These birational maps are actually pseudo-isomorphisms 
(\cite{_Huybrechts:cone_}, Proposition 4.7, 
 \cite{_Boucksom_}, Proposition 4.4), so that the second cohomologies
of $M$ and $M'$ are naturally identified. We shall sometimes say that
the  birational K\"ahler cone is the union of the K\"ahler cones of
birational models, omitting to mention pullbacks.

\hfill

\remark All this is a slight abuse of language: the birational K\"ahler cone
is not what one would normally call a cone (but its closure, the moving cone,
or birational nef cone, is). 


\hfill

\definition
{\bf An exceptional prime divisor} on a hyperk\"ahler manifold
is a prime divisor with negative square (with respect to the BBF form).

\hfill




The birational K\"ahler (or, rather, nef) cone is characterized in terms of prime 
exceptional divisors in the same way as the K\"ahler cone is characterized
in terms of rational curves.

\hfill

\theorem\label{birkahl-primeex} (Huybrechts, \cite{_Huybrechts:cone_}, Prop. 4.2)
Let $\eta\in \Pos(M)$ be an element of a positive cone
on a hyperk\"ahler manifold. Then $\eta$ is birationally nef
if and only if $q(\eta, E)\geq 0$ for any exceptional divisor
$E$. \endproof

\hfill

\remark In other words,  the faces of birational K\"ahler cone 
are dual to the classes of exceptional divisors.

\hfill

\theorem \label{_Markman_reflections_Theorem_} (Markman)
For each prime exceptional divisor $E$ on a hyperk\"ahler manifold,
there exists a reflection $r_E\in O(H^2(M,\Z))$ 
in the monodromy group, fixing $E^\bot$.

\hfill

{\bf Proof:} In the projective case, this is 
\cite[Theorem 1.1]{_Markman:reflections_}. The non-projective case
reduces to this by deformation theory. Indeed, as we shall show in the next section
(\ref{divisor-invariant}),
a prime exceptional divisor deforms locally as long as its
cohomology class stays of type $(1,1)$ (and its small deformations are
obviously prime exceptional again): in fact this also has already been 
remarked by
Markman in \cite{_Markman:reflections_}), but he restricted himself to the projective case.
Any hyperk\"ahler manifold $(M, I)$ has a small deformation $(M,I')$ which is projective.
The divisor $E$ deforms to $E'$ on $(M,I')$, so one obtains a reflexion 
$r_E'=r_E\in O(H^2(M,\Z))$ in the monodromy group (the monodromy being of topological 
nature, it is deformation-invariant in a tautological way).  

\endproof

\hfill

\definition
Such a reflection is called {\bf a divisorial reflection}.

\hfill

\definition
{\bf An exceptional chamber}, or {\bf a Weyl chamber} 
on a hyperk\"ahler manifold  is a connected component
of $\Pos(M)\backslash E^\bot$, where $E^{\bot}$ is a union of
all planes $e^\bot$ orthogonal to all prime exceptional divisors $e$. 

\hfill

\remark
An exceptional chamber is a fundamental domain of a group generated
by divisorial reflections. The birational K\"ahler cone is a dense open
subset of one of the exceptional chambers, which we shall call
{\bf birational K\"ahler chamber}.

\hfill

\theorem\label{_bir_mono_Theorem_} (Markman)\\
Let $(M,I)$ be a hyperk\"ahler manifold, and $\Mon^{\Hdg}(M,I)$ a
subgroup of the monodromy group fixing the Hodge
decomposition on $(M,I)$.
Then the image of $\Bir(M,I)$ in the orthogonal group of the Picard lattice 
is the group of all $\gamma\in \Mon^{\Hdg}(M,I)$ preserving the 
birational K\"ahler chamber $\Kah_B$.

\hfill

{\bf Proof:} Please see \cite[Lemma 5.11 (6)]{_Markman:survey_}.
The proof of \ref{_bir_mono_Theorem_} is similar to the proof of
\ref{_K3_auto_Theorem_}. 
\endproof

\hfill

Now we can generalize the birational Morrison-Kawamata cone
conjecture proved by Markman for projective hyperk\"ahler manifolds 
(\cite[Theorem 6.25]{_Markman:survey_}).

\hfill

\theorem\label{morkaw-birational}
Let $(M,I)$ be a hyperk\"ahler manifold, and $\Bir(M,I)$
the group of birational automorphisms of $(M,I)$.
Then $\Bir(M,I)$ acts on the set of faces of the
birational nef cone with finite number of orbits.

\hfill

{\bf Proof. Step 1:} Let $\delta$ be the discriminant of a lattice $H^2(M,\Z)$,
and $E$ an exceptional divisor. Then $|E^2|\leq 2\delta$.
Indeed, otherwise the reflection $x \arrow x - 2\frac{q(x,E)}{q(E,E)}E$
would not be integral.

{\bf Step 2:}  The group of isometries
of a lattice $\Lambda$ acts with finitely many orbits on the
set $\{l \in \Lambda\ \ |\ \ l^2 = x\}$
for any given $x$ (\ref{_Kneser_orbits_Theorem_}), and the monodromy group
is of finite index in the isometry group of the lattice $H^2(M,\Z)$.
Therefore,  $\Mon^{\Hdg}(M,I)$ is of finite index in the isometry group
of the Picard lattice and so all classes of exceptional divisors belong to
finitely many orbits of $\Mon^{\Hdg}(M,I)$.

{\bf Step 3:} We repeat the argument of \ref{_cone_for_K3_Corollary_}
to show that $\Mon^{\Hdg}(M,I)$ acts with finite number of orbits on the
set of faces of all exceptional chambers. It suffices to show that it
acts with finitely many orbits on the set of full flags
$P_{s-1}\supset P_{s-2}\supset \dots \supset P_1$ (notations as in 
the proof of \ref{_cone_for_K3_Corollary_}) formed by
intersections of orthogonal hyperplanes to the classes of exceptional
divisors. Using the fact that the squares of those classes are bounded
in absolute value by $C=2\delta$, we show that $P_{s-2}$ is described inside
$P_{s-1}$ as $y_1^\bot$ where $y_1$ is integral and $|y_1^2|\leq C^3$,
$P_{s-3}$ is defined in
$P_{s-2}$ as $y_2^\bot$ where $y_2$ is integral and $|y_1^2|\leq C^9$,
and so on. We deduce by \ref{_Kneser_orbits_Theorem_} that $O(\Lambda)$ acts 
with finitely many orbits
on the set of such full flags.

{\bf Step 4:} Thus $\Mon^{\Hdg}(M,I)$
acts with finite number of orbits on the
set of faces of all exceptional chambers.
For each pair of faces $F, F'$ of a birational K\"ahler
cone and $w\in O(\Lambda)$ mapping $F$ to $F'$, $w$ maps $\Kah_B$ to itself
or to an adjoint Weyl chamber $K'$. Then $K'=r(K)$, where
$r$ is the reflection fixing $F'$.
In the first case, $w\in \Aut(M)$.
In the second case, $rw$ maps $F$ to $F'$ and maps $\Kah_B$ to itself,
hence $rw\in \Aut(M)$. Therefore, there are as many
orbits of $\Mon^{\Hdg}(M,I)$ on the faces of exceptional
chambers as there are orbits of $\Bir(M,I)$.
\endproof

\hfill







\section{Deformation spaces of rational curves}
\label{_Defo_rac_Section_}


By a rational curve on a manifold $X$, we mean a curve $C\subset X$
such that its normalization is $\P^1$; in other words, the image of a
generically injective map $f:\P^1\to X$. Let $\Mor(\P^1, X)$ denote the
parameter space for such maps. Then by deformation theory 
(see \cite{_Kollar:curves_})
one has 
$$\dim_{[f]}(\Mor(\P^1, X))\geq \chi(f^*(TX))=-K_XC+\dim(X),$$
so that if $H$ denotes the space of deformations of $C$ in $X$, then
$$\dim(H)\geq -K_XC+\dim(X)-3.$$ 
A rational curve on an $m$-fold with trivial canonical class must therefore
move in a family of dimension at least $m-3$.

The following observation due to Z. Ran states that
on holomorphic symplectic manifolds, this estimate can be slightly improved.

\hfill

\theorem\label{Ran_ratcurves}
Let $M$ be a hyperk\"ahler manifold of dimension $2n$. Then any rational curve $C\subset M$ deforms 
in a family of dimension at least $2n-2$.

\hfill

{\bf Proof:} See \cite{Ran}, Corollary 5.2. Alternatively (we thank
Eyal Markman for this argument), one may notice that an extra parameter
is due to the existence of the twistor space $\Tw(M)$. This is a complex
manifold of dimension $n+1$, fibered over $\P^1$ in such a way that $M$
is one of the fibers and the other fibers correspond to the other complex 
structures coming from the hyperk\"ahler data on $M$. The map $f$, seen as a map
from $\P^1$ to $\Tw(M)$, deforms in a family of dimension at least $n+1$.
But all deformations have image contained in $M$ since the neighbouring fibers
contain no curves at all by \cite{_Verbitsky:trianal_} 
and the rational curves with
dominant projection to $\P^1$ belong to a different cohomology class.
\endproof

\hfill

Before making the following observation, let us recall a few definitions.

\hfill  

\definition
A complex analytic subvariety $Z$ of a holomorphically
symplectic manifold $(M, \Omega)$ is called {\bf
  holomorphic Lagrangian} if $\Omega\restrict Z=0$ and
$\dim_\C Z = \frac 1 2 \dim_\C M$, and {\bf isotropic} if
$\Omega\restrict Z=0$ (since $\Omega$ is non-degenerate, this implies $\dim_\C Z \leq \frac 1 2 \dim_\C
M$). It is called {\bf coisotropic} if $\Omega$ has rank
$\frac 1 2 \dim_\C M -\codim_\C Z$ on $TZ$ in all smooth
points of $Z$, which is the  minimal possible rank for a
$2n-p$-dimensional subspace in a $2n$-dimensional
symplectic space.

\hfill


Let now $Z$ be a compact K\"ahler manifold covered by rational curves.
By \cite{Campana-quotient}, there is an almost holomorphic\footnote{An almost 
holomorphic fibration is a rational map $\pi:\; X \arrow Y$,
holomorphic on $X_0\subset X$, and with $\pi^{-1}(y)\cap X_0$
non-empty for all $y\in Y$.} fibration
 $R: Z\dasharrow Q$,
called the {\bf rational quotient} of $Z$, such that its fiber through a sufficiently 
general point $x$ consists of all $y$ which can be joined from $x$ by a chain of rational curves.
That is, the general fiber of $R$ is {\bf rationally connected}. It is well-known that rationally
connected manifolds
do not carry any holomorphic forms, and it follows 
from \cite{Campana-quotient} that they are projective.

\hfill

\remark\label{GHS}
By a theorem due to Graber, Harris and Starr in the
projective setting, the base $Q$ of the rational 
quotient is not uniruled; this is also true in the compact
K\"ahler case (the reason being that the total space of a family
of rationally connected varieties over a curve is automatically
algebraic, so that the arguments of Graber-Harris-Starr apply), 
but we shall not need this.

\hfill

\theorem\label{_coiso_curves_Theorem_}
Let $M$ be a hyperk\"ahler manifold, $C\subset M$ a 
rational curve, and $Z\subset M$ be an irreducible component of the locus 
covered by the
deformations of $C$ in $M$. Then $Z$ is a coisotropic subvariety
of $M$. The fibers of the rational quotient of the desingularization of $Z$ have dimension
equal to the codimension of $Z$ in $M$.

\hfill


{\bf Proof:} We want to prove that at a general point of $Z$, the 
restriction of the symplectic form $\Omega$ to $TZ$ has kernel of
dimension equal to $k=\codim(Z)$. Recall that this kernel cannot be of dimension greater than $k$ by nondegeneracy of $\Omega$, and the equality means
that  $Z$ is coisotropic.

Let $h:Z'\to Z$ be the desingularization. The manifold $Z'$ is covered
by rational curves which are lifts of deformations of $C$. By dimension
count, through a general point $p$
of  $Z'$ there is a family of such curves of dimension at least $k-1$. In fact any rational 
curve through $p$ deforms in a family of dimension at least $k-1$, as seen by projecting it
to $M$ and applying \ref{Ran_ratcurves}. Take a minimal rational curve $C'$ through $p$,
then its deformations through $p$ cover a subvariety of dimension at least $k$.
This is because
by bend-and-break, there is only a finite number of minimal rational curves through two
general points (notice that bend-and-break applies here since all rational
curves through $p$ are in the fiber of the rational quotient, and this fiber
is projective). This means that the fibers of the rational quotient fibration
$R: Z'\dasharrow Q'$ are at least $k$-dimensional. But any holomorphic
form on $Z'$ is a pull-back of a holomorphic form on $Q'$ (since the 
rationally connected varieties do not carry any holomorphic
$m$-forms for $m>0$). So the tangent space to the
fiber of $R$ through $p$ is in the kernel of $h^*(\Omega)$, and thus the kernel of 
$\Omega|_{TZ}$ at $h(p)$ is of dimension $k$, and the same is true for the fiber of $R$ through $p$.
\endproof

\hfill

In the above argument, we have called a curve $C\subset Z$
{\bf minimal} if it is a rational curve of minimal degree
(say, with respect to the K\"ahler form
from the hyperk\"ahler structure) through a
general point of $Z$.

\hfill

\corollary\label{div-min}
If deformations of a rational curve $C$ in $M$ cover a divisor $Z$, then $C$ is a minimal rational curve 
in this divisor.

\hfill

{\bf Proof:} Indeed by 
\ref{_coiso_curves_Theorem_}, there is no other rational curve
through a general point of $Z$.
\endproof

\hfill

\corollary\label{min-dim}
Minimal rational curves on holomorphic symplectic manifolds deform in a 
$2n-2$-parameter family. 

\hfill

{\bf Proof:} This is obvious from the proof of \ref{_coiso_curves_Theorem_}.
Indeed, if the deformations of a minimal $C$ cover a subvariety $Z$ of 
codimension $k$, then these deformations form a family of dimension
$k-1+dim(Z)-1=k-1+2n-k-1=2n-2$, since there is a $k-1$-parametric family
of them through the general point of $Z$.
\endproof

\hfill

\remark\label{contractibles}
These results have well-known analogues (which are consequences of the work by Wierzba and Kaledin) 
in the case when $Z$ is 
contractible, that is, one has a birational morphism $\pi:M\to Y$ whose exceptional set is $Z$. 
The image of $Z$ by $\pi$ replaces 
the base of the rational quotient. In the case when moreover $Z$ is a 
divisor, the 
fibers of $\pi$ are one-dimensional and therefore are trees of smooth rational 
curves by Grauert-Riemenschneider theorem. In general, it is not obvious
whether the minimal rational curves are smooth. In one important case, though,
they are smooth: namely, when $Z$ is a negative-square divisor on $M$. The 
reason
is that one can deform the pair $(M,Z)$ so that $M$ becomes projective
(see below), and
then use results by Druel from \cite{druel-contract} to reduce to the
contractible case.

In this case, the inverse image of the tangent bundle $T_M$ on the
normalization of a minimal rational curve $C$ splits as follows:
$$f^*T_M={\cal O}(-2)\oplus{\cal O}(2)\oplus{\cal O}\oplus\dots\oplus{\cal O}$$ This is because it should be isomorphic to its dual (by the pull-back of
the symplectic form), have at most one negative summand, and contain 
$T_{\P^1}$ as a subbundle. This also remains true when the normalization map 
of $C$
is an immersion. In general, there are problems related to the singularities
(for instance, as soon as there are cusps, $T_{\P^1}$ is only a subsheaf and 
not a subbundle
of $f^*T_M$) and it is not obvious how to avoid them in order to get a
similar splitting, which one would like to be 
$f^*T_M={\cal O}(-2)\oplus{\cal O}(2)\oplus{\cal O}^{\oplus k}
\oplus({\cal O}(-1)\oplus{\cal O}(1))^{\oplus l}$. 

\hfill

\corollary\label{min-deforms} If $C$ is minimal, any small deformation $M_t$ 
of $M=M_0$ such that
the dual class $z$ of $C$ stays of type $(1,1)$ on $M_t$, contains a 
deformation of $C$.

\hfill

{\bf Proof:}
Consider the universal family ${\cal M}\to \Def(M)$ of small deformations of $M=M_0$.
The $t\in \Def(M)$ such that $z$ is of type $(1,1)$ on $M_t$ form a smooth 
hypersurface in 
$\Def(M)$ (if one identifies the tangent space to $\Def(M)$ at $0$ with 
$H^1(M, T_M)\cong
H^1(M, \Omega^1_M)$, then the tangent space to this hypersurface is the
hyperplane orthogonal to $z$).
Now let ${\cal C}\to \Def(M)$ be the family of deformations of 
$C$ in ${\cal M}$: the image of ${\cal C}$ in $\Def(M)$ is a subvariety
by Grauert's proper mapping theorem, and
it suffices to prove
that it is a hypersurface. This is a simple
dimension count. Indeed one obtains from Riemann-Roch theorem, as in 
the proof of \ref{Ran_ratcurves}, that $C$ deforms in a family of dimension
at least $2n-3+\dim(\Def(M))$. Since the deformations of $C$ inside any $M_t$
form a family of dimension $2n-2$ (when nonempty), the conclusion follows.
\endproof

\hfill

\corollary\label{min-deformsglobally} If $C$ is minimal, any
deformation of $M=M_0$ such that the corresponding homology class
remains of type $(1,1)$, has a birational model containing a rational
curve in that homology class.

\hfill

{\bf Proof:} This is the same argument as in the previous corollary,
but we have to consider the universal family over $\Teich_z(M)^0$,
where $\Teich(M)^0$ is the connected component of $\Teich(M)$ containing 
the parameter point for our complex manifold $M_0$, and $\Teich_z(M)^0$
is the part of it where $z$ remains of type $(1,1)$.
Birational models appear since $\Teich_z(M)$ is not Hausdorff, so that
a subvariety of $\Teich_z(M)^0$ of maximal dimension is not necessarily
equal to $\Teich_z(M)^0$; on the other hand it is known by the work of
Huybrechts that unseparable points of $\Teich_z(M)$ correspond to 
birational complex manifolds.
\endproof

\hfill


The deformations of a minimal $C$ as above are obviously minimal on the
neighbouring fibers (indeed, new effective classes only appear on closed subsets of
the parameter space). However,
globally a limit of minimal curves does not have 
to be minimal, at least apriori; it can also become reducible (and in fact 
does so
in many examples). The deformation theory of non-minimal 
curves is not as nice
as described in \ref{min-deforms}. For instance, in an 
example from \cite{bayer-hass-tschi}, credited by the authors to Claire 
Voisin, there is a holomorphic symplectic fourfold $X$ containing a lagrangian
quadric $Q$, and rational curves of type $(1,1)$ on $Q$ deform only together
with $Q$, that is over a codimension-two subspace of the base space,
as follows from \cite{voisin-lagrang}. However this problem disappears,
at least locally, 
in the case when the dual class $z$ is {\bf negative divisorial},
that is, when $q(z,z)<0$ and $z$ is $\Q$-effective.

The following theorem is due to Markman in the projective case.

\hfill

\theorem\label{divisor-invariant} Let $z$ be a negative $(1,1)$-class
on $M$ which is represented by an irreducible divisor $D$. Then $z$ is 
effective on any deformation $M_t$ where it stays of type $(1,1)$, and 
represented by an irreducible divisor on an open part of this parameter 
space.

\hfill

{\bf Proof:} By \cite{_Boucksom_}, $D$ is uniruled. By 
\ref{_coiso_curves_Theorem_}, there is only one dominating family of rational
curves on $D$. Take a general member $C$ in this family. We have seen that
it deforms on all neighbouring $M_t$, yielding rational curves $C_t$. In 
particular its dual class is proportional (with some positive coefficient) 
to that of $D$, since both stay
of type $(1,1)$ on the same small deformations of $M$; therefore
$CD<0$ and all deformations of $C$ in $M$ stay inside $D$. We claim that the deformations of 
$C_t$ on $M_t$ also cover divisors on $M_t$, which have then to be 
deformations of $D$. Indeed, let $Z_t$ be the subvariety of $M_t$ covered
by the curves $C_t$. Suppose it is not a divisor for all $t$. Since the
dimension jumps over closed subsets of the parameter space, one has an open
subset $U$ of such $t$. By \ref{_coiso_curves_Theorem_},
through a general point of $Z_t, t\in U,$ there is a 
positive-dimensional family of curves $C_t$. Now consider a parameter space
for the following triples over $T$: $\{(C_t, x, C'_t): x\in C_t\cap C'_t\}$.
The dimension of its fibers over $t\in U$ is strictly greater than the
dimension of the central fiber, and this is a contradiction. 

Therefore $D$ deforms everywhere locally. Now consider the universal 
family ${\cal M}$  over $\Teich_z$, and the ``universal divisor'' 
${\cal D}$ in this family.
The image of ${\cal D}$ is a subvariety in $\Teich_z$, of
the same dimension as $\Teich_z$. This means that any deformation 
of $M$ preserving the type $(1,1)$ of $z$ has a birational model such that
$z$ is effective. As birational hyperk\"ahler manifolds differ only
in codimension two, the theorem follows.  
\endproof

\hfill

\corollary\label{divisor-pic1} A negative class $z\in H^2(M,\Z)$ is 
effective 
or not
simultaneously in all complex structures where it is of type $(1,1)$
and generates the Picard group over $\Q$.

\hfill

{\bf Proof:} If the only integral $(1,1)$-classes are rational multiples 
of $z$, one cannot have more than one effective divisor on $M$: indeed if 
$D$ and $D'$ are prime effective divisors and $q(D,D')<0$,
then $D=D'$ (if not, remark that by definition $q(D,D')$ is obtained
by integration of a positive form over $D\cap D'$). In particular, every
effective divisor is irreducible. Now apply the previous theorem.
\endproof

\hfill

When one wants to deform curves rather than divisors, the notion of a
minimal curve that we have used above is not quite natural since it 
depends, apriori, on the subvariety $Z\subset M$ covered by deformations
of the curve, and not only on the complex manifold $M$ itself. In the
projective setting, one typically considers rational curves generating 
an extremal ray of the Mori cone: the class of such a curve is {\bf extremal}
in the sense of the Introduction. In the non-projective setting, the
following notion of extremality looks better-behaved.

\hfill

\definition\label{minimal} Let $X$ be a K\"ahler manifold and $z$ an integral homology 
class of type
$(1,1)$. We say that $z$ is {\bf minimal} if the intersection of $z^\bot$ with
the boundary of the K\"ahler cone contains an open subset of $z^\bot$.

\hfill

An example due to Markman, \cite{_Markman:reflections_},
shows that a limit of minimal
(or extremal) curves does not have to be minimal
(extremal).

\hfill

\example (\cite{_Markman:reflections_}, Example 5.3) Let $\bar X_0$ be an intersection of a quadric
and a cubic in $\P^4$ with one double point. The resolution 
$p:X_0\to \bar X_0$ is a
$K3$ surface, and $p^*H+2E$, where $H$ is a hyperplane section and $E$ is
the exceptional curve, is not minimal (or extremal). Now deform $X_0$ to
a smooth non-projective K3 surface $X_t$ in such a way that only the multiples 
of $p^*H+2E$ survive in 
$H^{1,1}$. Since $p^*H+2E$ is a $(-2)$-class, it is easy to show that those are effective. They are obviously
extremal (and minimal, too).

\hfill

In what follows, we shall define the {\bf MBM-classes} as classes which are
minimal modulo monodromy and birational equivalence. We shall show that
this notion is deformation-invariant, and that if a class is MBM, then
the intersection of its orthogonal hyperplane and the positive cone is 
a union of faces of K\"ahler-Weyl chambers 
(\ref{_KW_chambers_MBM_Theorem_}, \ref{MBM-invar-thm}). It is an interesting question whether the notions of minimality 
and extremality 
are equivalent
for hyperk\"ahler manifolds. Minimal classes are extremal, and it follows from our proofs 
that minimal classes are effective up to a rational multiple and represented by rational curves, as 
it is the case
of extremal rays in the Mori theory. But extremal classes do not, apriori, have to be minimal,
since extremal rays of the cone of curves could accumulate. One can conjecture that in real 
life it never happens (as this is true in the projective case by \cite{HT-GAFA}).


\section{Twistor lines in the Teichm\"uller space}


\subsection{Twistor lines and 3-dimensional planes in $H^2(M,\R)$.}
\label{_twi_lines_Subsection_}

Recall that any hyperk\"ahler structure $(M,I,J,K,g)$
defines a triple of K\"ahler forms 
$\omega_I, \omega_J, \omega_K\in \Lambda^2(M)$
(Subsection \ref{_hk_twi_Subsection_}). A hyperk\"ahler structure
on a simple hyperk\"ahler manifold is determined by a 
complex structure and a K\"ahler class (\ref{_Calabi-Yau_Theorem_}).

\hfill

\definition
Each hyperk\"ahler structure induces a family $S\subset \Teich$ of deformations
of complex structures parametrized by $\C P^1$ (Subsection
\ref{_hk_twi_Subsection_}). The curve $S$ is called
{\bf the twistor line} associated with the
hyperk\"ahler structure $(M,I,J,K,g)$.

\hfill

We identify the period space $\Perspace$ with the Grassmannian
of positive oriented 2-planes in $H^2(M,\R)$ 
(\ref{_period_Grassmann_Proposition_}).
For any point $l\in S$ on the twistor line, the
corresponding 2-dimensional space $\Per(l)\in \Perspace$
is a 2-plane in the 3-dimensional space 
$\langle \omega_I, \omega_J, \omega_K\rangle$.
Therefore, $\langle \omega_I, \omega_J, \omega_K\rangle$
is determined by the twistor line $S$ uniquely, as the linear
span of the planes in $\Per(S)$.

We call two hyperk\"ahler structures {\bf equivalent}
if one can be obtained from the other by a homothety and a 
quaternionic reparametrization:
\[ (M,I,J,K,g)\sim (M,hIh^{-1},hJh^{-1},hKh^{-1},\lambda g),\]
for $h\in {\Bbb H}^*$, $\lambda\in \R^{>0}$. 
Clearly, equivalent hyperk\"ahler structures produce the same
twistor lines in $\Teich$. However, a 
hyperk\"ahler structure is determined by a complex
structure, which yields a 2-dimensional subspace $\Per(I)=\langle \omega_J,\omega_K\rangle$
in $H^2(M,\R)$, and a K\"ahler structure $\omega_I$,
as in \ref{_Calabi-Yau_Theorem_}. The form $\omega_I$ can be
reconstructed up to a constant from the 3-dimensional space
$\langle \omega_I,\omega_J,\omega_K\rangle$ and the plane
$\Per(I)$. This proves the following result, which is
essentially a form of Calabi-Yau theorem.

\hfill

\claim
Let $(M,I,J,K,g)$ be a hyperk\"ahler structure on
a compact manifold,
and $S\subset \Teich$ the corresponding twistor line. Then $S$
is sufficient to recover the equivalence class
of $(M,I,J,K,g)$. 
\endproof

\hfill

\proposition\label{_3-plane_Teich_z_Proposition_}
Let $S\subset \Teich$ be a twistor line,
$W\subset H^2(M,\R)$ be the corresponding 3-dimensional
plane $W:=\langle \omega_I, \omega_J, \omega_K\rangle$,
and $z\in H^2(M,\R)$ a non-zero real cohomology class. 
Then $S$ lies in $\Teich_z$ if and only if $W\bot z$.

\hfill

{\bf Proof:} By definition, $\Teich_z$ is the set
of all $I\in \Teich$ such that the 2-plane $\Per(I)$
is orthogonal to $z$. However, any point of $\Per(S)$
lies in $W$, hence all points in $S$ belong to $\Teich_z$ 
whenever $W\bot z$. Conversely, the planes corresponding to 
points in $\Per(S)$ generate $W$, so that $W$ is orthogonal
to $z$ if all those planes are.
\endproof

\subsection{Twistor lines in $\Teich_z$.}

\remark\label{_Teich_z_twistor_curve_Remark_}
Let $z\in H^{1,1}(M,I)$ be a non-zero cohomology class on 
a hyperk\"ahler manifold $(M,I)$, $\Teich^I$ the connected
component of the Teichm\"uller space containing $I$, 
and $\Teich_z$ the set of all $J\in \Teich^I$ such that
$z$ is of type $(1,1)$ on $(M,J)$. Given a K\"ahler
form $\omega$ on $(M,I)$ such that  $q(\omega,z)=0$,
consider the corresponding hyperk\"ahler structure
$(M,I,J,K)$, and let $S\subset \Teich$ be the
corresponding twistor line. By \ref{_3-plane_Teich_z_Proposition_}, 
$S$ lies in $\Teich_z$. 

In other words, there is a twistor line on $\Teich_z$  through the point $I$
if and only if $z$ is orthogonal to a K\"ahler form on $(M,I)$.

\hfill

\definition
In the assumptions of \ref{_Teich_z_twistor_curve_Remark_},
let $W_S:=\langle \omega_I, \omega_J,\omega_K\rangle$
be a 3-dimensional plane associated with the hyperk\"ahler
structure $(M,I,J,K)$, and $S\subset \Teich_z$ 
the corresponding twistor line. The twistor line $S$
is called {\bf $z$-GHK} ($z$-general hyperk\"ahler) if
$W^\bot\cap H^2(M,\Q)=\langle z\rangle$.

\hfill

The utility of $z$-GHK lines is that they can be lifted
from the period space to $\Teich_z$ uniquely, if $\Teich_z$
contains a twistor line (\ref{_Per_on_hk_surge_z_GHK_Corollary_}).
For other twistor lines, such a lift, even if it exists, is not necessarily
unique because of non-separable points, but 
a $z$-GHK line has to pass through a Hausdorff 
point of $\Teich_z$.

\hfill

\claim\label{admits-GHK} 
Let $z\in H^{1,1}(M,I)$ be a non-zero cohomology class on 
a hyperk\"ahler manifold $(M,I)$. Assume that $(M,I)$
admits a  K\"ahler form $\omega$ such that  $q(\omega,z)=0$. 
Then $(M,I)$ admits a $z$-GHK line.

\hfill

{\bf Proof:} Let $z^\bot\subset H^{1,1}_I(M,\R)$ be 
the orthogonal complement of $z$. The set of K\"ahler
classes is non-empty and open in $z^\bot$.
For each integer vector $z_1\in H^2(M,\Z)$,
non-collinear with $z$, the orthogonal complement
$z_1^\bot$ intersects with $z^\bot$ transversally.
Removing all the hyperplanes $z^\bot \cap z_1^\bot$
from $z^\bot$, we obtain a dense set which
containes a K\"ahler form $\omega_1$, such that
$\omega_1^\bot$ contains none of the integer 
vectors $z_1$. The hyperk\"ahler structure
associated with $I$ and $\omega_1$
satisfies $W^\bot\cap H^2(M,\Q)=\langle z\rangle$,
because $W\ni \omega_1$, and 
$\omega_1^\bot\cap H^2(M,\Q)=\langle z\rangle$.
\endproof

\hfill

The following result is easy to prove,
however, we give a reference to \cite{_V:Torelli_},
where the proof is spelled out in a situation which
is almost the same as ours.

\hfill

\claim\label{_GHK_lines_and_Pic_Claim_}
Let $S\subset \Teich_z$ be a twistor line,
associated with a 3-\-di\-men\-sio\-nal subspace 
$W\subset H^2(M, \R)$. 
Then the following assumptions are equivalent.
\begin{description}
\item[(i)] $S$ is a $z$-GHK line.
\item[(ii)] For all $l\in S$, except a countable number,
the space $H^{1,1}(M_l,\Q)$ is generated by $z$, where
$M_l$ is the manifold corresponding to $l$.
\item[(iii)] For some $w\in W$, the space of rational
classes in the orthogonal complement 
$w^\bot\subset H^2(M, \R)$ is generated by $z$.
\end{description}
{\bf Proof:} \cite[Claim 5.4]{_V:Torelli_}. \endproof

\hfill

Given a $z$-GHK line $S$, we may chose a point $l\in S$
where $H^{1,1}(M_l,\Q)=\langle z\rangle$. For this point,
the K\"ahler cone coincides with the positive cone by
\ref{Boucksom-cone}. Indeed, $M_l$ contains no curves, as the class $z$,
being orthogonal to a K\"ahler form, cannot be effective.
  Chosing a different K\"ahler form in 
$\Pos(M_l)\cap z^\bot$, we obtain a different twistor
line in $\Teich_z$, intersecting $S$. It was shown that
such a procedure can be used to connect any two points
of $\Teich_z$, up to non-separatedness issues (in \cite{_V:Torelli_} it was proven
for the whole $\Teich$, but for $\Teich_z$ the argument is
literally the same).

\hfill

\proposition\label{_4_GHK_lines_Proposition_}
Let $x, y \in \Teich_z$. Suppose that $\Teich_z$ contains
a twistor line. Then a point $\tilde x$ non-separable from
$x$ can be connected
to a point $\tilde y$ non-separable from
$y$ by a sequence of at most 5 sequentially intersecting
$z$-GHK lines.

\hfill

{\bf Proof:} From \cite[Proposition 5.8]{_V:Torelli_} 
it follows that any two points $a=\Per(x), b=\Per(y)$ in the period domain
$\Per(\Teich_z)$ can be connected by at most 5 sequentially intersecting
$z$-GHK lines. Moreover, the intersection points can be
chosen generic, in particular, separable. Lifting
these $z$-GHK lines to $\Teich_z$, we connect
a point $\tilde x$ in $\Per^{-1}(x)$ to a point $\tilde y$
in $\Per^{-1}(y)$ by a sequence of at most 5 sequentially intersecting
$z$-GHK lines. However, these preimages are non-separable
from $x$, $y$ by the global Torelli theorem (\ref{_glo_Torelli_Theorem_}).
\endproof

\hfill

\remark Clearly, if our twistor line already passes through $x$, we can connect
$x$ itself to a $\tilde y$ non-separable from $y$.

\hfill

\subsection{Twistor lines and MBM classes.}

As we have already mentioned in the Introduction, the Bogomolov-Beauville-Fujiki form identifies 
the rational cohomology and homology of a 
hyperk\"ahler manifold $M$, inducing an injective map
$q: H^2(M,\Z)\to H_2(M,\Z)$. Since $q$
is not necessarily unimodular, this map is 
not an isomorphism over $\Z$. However, it is
an isomorphism over $\Q$, and we shall identify
$H^2(M,\Q)$ with $H_2(M,\Q)$ without further comment
(in particular, we shall often view the classes in $H_2(M,\Z)$ as
rational classes in $H^2$). The only exception to this rule is the
notion of effectiveness: an effective homology class is a class of a curve,
whereas an effective cohomology class is a class of an effective divisor.

Recall that a negative integer homology class
$\eta\in H_2(M,\Z)$ is called an MBM class if 
some image of $\eta$ under monodromy contains in its
orthogonal a face
of the K\"ahler cone of some birational hyperk\"ahler model $M'$ 
(which is necessarily pseudo-isomorphic to $M$).  To study the MBM classes
it is convenient to work with hyperk\"ahler manifolds
which satisfy $\rk \Pic(M)=1$, with $\Pic(M)\subset H^2(M,\Z)$
generated over $\Q$ by $\eta$. Here $\eta$ is considered as an element of
$H^2(M,\Q)$, identified with $H_2(M,\Q)$ as above.
In this case $M$ is clearly non-algebraic since $\eta$ is negative.

\hfill

\theorem\label{_MBM_posi_theorem_}
Let $(M,I)$ be a hyperk\"ahler manifold, such that 
$\Pic(M,I)=\langle z \rangle$, where
$z\in H_{1,1}(M,I)$ is a non-zero negative class. Then the following 
statements are equivalent.
\begin{description}

\item[(i)] The class $\pm z$ is $\Q$-effective.
\item[(ii)] The class $\pm z$ is extremal.
\item[(iii)] The K\"ahler cone of $(M,I)$ is not equal
to its positive cone.
\item[(iv)] The class $z$ is minimal, in the sense of \ref{minimal}.
\item[(v)] The class $z$ is MBM.
\end{description}

{\bf Proof:} First of all, the equivalence of (i) and (ii) is clear since the cone of curves
on $M$ consists of multiples of $\pm z$ if it is $\Q$-effective and is empty otherwise.

The equivalence of (i) and (iii) is 
proven as follows: first of all, for signature reasons
one can find a positive class $\alpha$ with $q(\alpha,z)=0$, and if $\pm z$ is $\Q$-effective
then $\alpha$ cannot be K\"ahler. The converse follows by Huybrechts-Boucksom's description
of the K\"ahler cone, see \ref{Boucksom-cone} (if $\pm z$ is not $\Q$-effective then there are no curves at all, in
particular no rational curves, and thus any positive class is K\"ahler).

Finally, let us show the equivalence of (iii), (iv) and (v). 

If $\Kah(M,I)\neq \Pos(M,I)$, then by Huybrechts-Boucksom 
the K\"ahler cone has faces supported on hyperplanes
of the form $x^\bot$ where $x$ is a class of a rational curve.
However, the only integral $(1,1)$-classes on $(M,I)$ are multiples of $z$,
so such a face can only be $z^\bot$. In this case $z$ must be minimal;
this proves (iii) $\Rightarrow$ (iv). The implication (iv) $\Rightarrow$ (v)
is a tautology. Finally, if $\Kah(M,I)=\Pos(M,I)$, there are no minimal classes
and no non-trivial birational models of $M$, so there cannot be any MBM classes
either, which proves (v) $\Rightarrow$ (iii). \endproof

\hfill





\hfill

In general, the notion of an MBM class is much more complicated,
but one can hope to study them by deforming to the case we have just
described. This is indeed what we are going to do, using 
\ref{_4_GHK_lines_Proposition_} as a key tool. The following useful proposition,
which allows one to identify generic positive 3-dimensional subspaces of 
$H^2(M,\R)$ and twistor lines, is an illustration of our method.

\hfill

\proposition\label{_Per_on_hk_surge_z_GHK_Corollary_}
Let $z\in H^2(M,\Q)$ be a non-zero vector, such that
$\Teich_z$ contains a twistor line, and $W\subset z^\bot$
a 3-dimensional positive subspace of $H^2(M,\R)$ which satisfies 
$W^\bot \cap H^2(M,\Q)\subset \langle z\rangle$.
Then there exists a $z$-GHK twistor line $S$ such that
the corresponding 3-dimensional space $W_S\subset H^2(M,\R)$
is equal to $W$.

\hfill

{\bf Proof:} Let $V\subset W$ be a 2-dimensional plane
which satisfies \[ V^\bot \cap H^2(M,\Q)\subset \langle z\rangle\]
(by \ref{_GHK_lines_and_Pic_Claim_} (ii) 
this is true for all planes except at most 
a countable set). Then $V=\Per(I)$, where $I\in
\Teich_z$, because the period map is surjective.
The space $H^{1,1}_\Q(M,I)=V^\bot\cap H^2(M,\Q)$
is generated by $z$. Since $\Teich_z$ contains a twistor line,
by \ref{_4_GHK_lines_Proposition_} there must be a twistor line 
through the point $I'\in \Teich_z$ which is inseparable from $I$.

Now on $I'$, $z$ cannot be effective since it is orthogonal to a K\"ahler form.
But this means that there are no curves on $(M,I')$, so the K\"ahler cone
of $(M,I')$ is equal to the positive cone, $\Teich_z$ is separated at $I'$,
and $I=I'$. Therefore the K\"ahler cone of $(M,I)$ coincides with its
positive cone, hence there exists a K\"ahler
form $\omega$ on $(M,I)$ such that $\langle V, \omega\rangle = W$.
The corresponding hyperk\"ahler structure gives
a 3-plane $W_S=\langle V, \omega\rangle = W$.
\endproof

\hfill

Our strategy in proving the deformation invariance of MBM property 
is as follows: we first show that the property of being
orthogonal to a K\"ahler form, modulo monodromy and birational
transformations, is deformation-invariant, and then show that the MBM classes
are exactly those which do not have this property.

\hfill

\theorem\label{MBM-invar-thm} Let $M$ be a simple hyperk\"ahler manifold, 
and $z\in H^2(M,\Q)$
a cohomology class. Suppose that 
$z$ is orthogonal to a K\"ahler form on $(M,I_0)$, where
$I_0\in \Teich_z$. Then:
\begin{description}
\item[(i)]  
For any $I\in \Teich_z$, there is some $I'$ non-separable from $I$ such that 
$z$ is orthogonal to a K\"ahler form on $(M,I')$.

\item[(ii)] On the manifold $(M,I)$, the class
$z$ is orthogonal to an element of a K\"ahler-Weyl chamber (that is, 
 $z$ is orthogonal to $\gamma(\alpha)$ for some $\gamma$ in the
monodromy group and $\alpha$ lying in the birational K\"ahler cone, i.e. 
K\"ahler on a birational model, see \ref{chamber}). 
\item[(iii)] The elements of K\"ahler-Weyl chambers are dense in the intersection of $z^\bot$
and the positive cone.
\end{description}

\hfill 

{\bf Proof:} (i) By  \ref{_GHK_lines_and_Pic_Claim_}, there is a $z$-GHK line through $I_0$. By 
\ref{_4_GHK_lines_Proposition_} ,
$I_0$ is connected to $I'$ by a sequence of $z$-GHK lines. Therefore on $(M,I')$, $z$ is 
orthogonal to a K\"ahler form $\omega'$.

(ii) The group $\Mon^{\Hdg}(M,I)$ acts transitively on the set of the Weyl chambers (see \cite{_Markman:survey_};
in fact this easily follows from \ref{_Markman_reflections_Theorem_}). Let $W(I)$ denote the fundamental
Weyl chamber, that is, the interior of the birational nef cone. Then there is an element 
$\gamma\in \Mon^{\Hdg}(M,I)$
such that $W(I)=\gamma W(I')$.

Consider $\tilde\gamma$ which is a lift of $\gamma$ in the mapping class group $\Gamma_I$, and the
complex structure $\tilde\gamma(I')$. This is again non-separable from $I$ and $I'$ and it carries
a K\"ahler form $\tilde\gamma(\omega')$, orthogonal to $\gamma(z)$. Its cohomology class is an
element of the birational K\"ahler cone of $I$. The inverse image of this class by $\gamma$ 
is orthogonal to $z$, q.e.d..

(iii) We claim that any positive $(1,1)$-class $\omega$ which is orthogonal to $z$ but not to any other rational
cohomology class is K\"ahler modulo monodromy and
birational transforms. Indeed, consider the twistor line 
$L\subset \Perspace$
corresponding to 
the 3-subspace generated by the period of $I_0$ and $\omega$. By 
\ref{_Per_on_hk_surge_z_GHK_Corollary_},
it lifts to $\Teich_z$ as a $z$-GHK line. 
This line does not have to pass through $I_0$, but it passes through 
some $I'_0$ nonseparable from $I_0$. On this $I'_0$,
$\omega$ is K\"ahler, and we conclude in the same way as 
in (ii) that on $I_0$ itself it is K\"ahler
modulo monodromy and birational equivalence. \endproof

\hfill

Notice that by definition, the property obtained in (iii) expresses exactly the fact that $z$ is not MBM.
We thus obtain the deformation-invariance of MBM classes.

\hfill

\corollary\label{MBM-invar-coro} A negative class $z$ is MBM or not simultaneously  
in all complex structures where it is of type $(1,1)$.
\endproof

\hfill

\corollary\label{MBM-effec} An MBM class is $\pm \Q$-effective and represented by a 
rational curve 
(up to a scalar) on a birational model of $(M,I)$.

\hfill

{\bf Proof:} We have seen that such is the case (even on $(M,I)$ itself) when $z$ generates the Picard group over $\Q$. The rest follows by
deformation invariance of MBM property and the results on deformations of minimal rational curves in section 4. \endproof

\hfill

Putting all we have done in this section together, we arrive at the following.

\hfill

\theorem\label{MBM-characterization}
Let $M$ be a simple hyperk\"ahler manifold, and $z\in H_2(M,\Q)$
a homology class. Then the following statements are equivalent.
\begin{description}
\item[(i)] The space $\Teich_z$ (the subset of all $I\in \Teich$ such that
$z$ lies in $H_{1,1}(M,I)$) contains a twistor  line.
\item[(ii)] For each $I\in \Teich_z$,
there exists $I'\in \Teich_z$ non-separable 
from $I$ which is contained in a twistor line.
\item[(iii)] For each $I\in \Teich_z$ with Picard number one, $\pm z$ is not effective;
\item[(iv)] For some $I\in \Teich_z$ with Picard number one, $\pm z$ is not effective;
\item[(v)] $z$ is not MBM.
\end{description}
Moreover in the items (iii) and (iv) one can replace ``$\pm z$ is not effective'' by
``$(M,I)$ contains no rational curves''.

\hfill

{\bf Proof:} The implication (i) $\Rightarrow$ (ii) follows from
\ref{_4_GHK_lines_Proposition_}. If (ii) holds, then there is $I'$ non-separable from $I$
such that $z$ is orthogonal to a K\"ahler form on $I'$, but, as we have already seen, 
one then has $\Kah(I')=\Pos(I')$ and $I=I'$, so $z$ cannot be effective on $I'$, hence (iii).
(iii) $\Rightarrow$ (iv) is obvious. To get (i) from (iv), notice that  $\Kah(I)=\Pos(I)$ so there
are K\"ahler classes orthogonal to $z$ on $(M,I)$.
This implies existence of twistor lines in $\Teich_z$.
The property (v) is equivalent to (iv) since 
these properties are deformation-invariant, and the
equivalence for complex structures with Picard number one has already been verified.
Finally, if $z$ is effective on $(M,I)$ with Picard number one, it is automatically represented by
a rational curve. Indeed otherwise there are no rational curves at all, so the K\"ahler cone
should be equal to the positive cone, but one easily finds a positive $(1,1)$-form orthogonal 
to $z$. \endproof

\section{Monodromy group and the K\"ahler cone}


In this section we prove the results on the K\"ahler cone stated in
the Introduction.

\subsection{Geometry of K\"ahler-Weyl chambers}

\definition\label{chamber}
Let $(M,I)$ be a hyperk\"ahler manifold, and $\Mon^{\Hdg}(M,I)$
the group of all monodromy elements preserving the Hodge decomposition
on $(M,I)$. A {\bf K\"ahler-Weyl chamber} of a hyperk\"ahler manifold
is the image of the K\"ahler cone of $M'$ under some 
$\gamma\in \Mon^{\Hdg}(M,I)$,
where $M'$ runs through the set of all birational models of $M$.

\hfill

\theorem\label{_KW_chambers_MBM_Theorem_}
Let $(M,I)$ be a hyperk\"ahler manifold, 
and $S\subset H_{1,1}(M,I)$ the set of all MBM classes in 
$H_{1,1}(M,I)$. Consider the corresponding set of hyperplanes
$S^\bot:=\{W=z^\bot\ \ |\ \ z\in S\}$ in $H^{1,1}(M,I)$.
Then the K\"ahler cone of $(M,I)$ is 
a connected component of $\Pos(M,I)\backslash S^\bot$,
where $\Pos(M,I)$ is a positive cone of $(M,I)$.
Moreover, the connected components of $\Pos(M,I)\backslash S^\bot$
are K\"ahler-Weyl chambers of $(M,I)$. 

\hfill



{\bf Proof}: First of all, none of the classes $v\in W,\ W\in S^\bot$
belong to interior of any of the K\"ahler-Weyl chambers:
indeed it follows from \ref{MBM-invar-thm} that if $v\in W$
is in the interior of a K\"ahler-Weyl chamber, then the classes
belonging to the interior of K\"ahler-Weyl chambers are dense 
in $W$ so $W\not\in S^\bot$. 

It remains to show that any connected component 
of $\Pos\backslash S^\bot$ is a K\"ahler-Weyl chamber.

Consider now a cohomology class $v\notin S^\bot$,
and let $W=\langle \Re\Omega,\Im\Omega, v\rangle$
be the corresponding 3-dimensional plane in $H^2(M,\R)$.
We say that $v$ is KW-generic if $W^\bot \cap H^{1,1}(M,\Q)$ is
at most 1-dimensional. Denote by ${\goth W}$ the set of 
non-generic $(1,1)$-classes. Clearly, ${\goth W}$ is a union of
codimension 2 hyperplanes, hence removing ${\goth W}$ from
$\Pos\backslash S^\bot$ does not affect the set of connected
components.

We are going to show that any 
$v\in \Pos(M,I)\backslash ({\goth W}\cap S^\bot)$
belongs to the interior of a K\"ahler-Weyl chamber. 
This would imply that the connected components of
$\Pos(M,I)\backslash  S^\bot$ are K\"ahler-Weyl chambers,
finishing the proof of  \ref{_KW_chambers_MBM_Theorem_}.

Since $v$ is KW-generic, and not in $S^\bot$, 
the space $W^\bot \cap H^{1,1}(M,\Q)$ contains no 
MBM-classes. This implies that $\Teich_z$ is covered
by twistor lines, where $z$ is a generator of
$W^\bot \cap H^{1,1}(M,\Q)$. By \ref{_Per_on_hk_surge_z_GHK_Corollary_},
there exists a hyperk\"ahler structure
$(I',J,K)$ such that $W$ is the corresponding
3-dimensional plane
$\langle \omega_{I'}, \omega_J, \omega_K\rangle$
and $\Per(I')=\Per(I)$. Therefore, $v$ is K\"ahler
on $(M,I')$. As we have already seen in the proof of \ref{MBM-invar-thm},
this implies that $v$ belongs to  the interior of a
K\"ahler-Weyl chamber on $(M,I)$. \endproof

\subsection[Morrison-Kawamata cone conjecture  and minimal curves]{Morrison-Kawamata cone conjecture \\and minimal curves}

Recall the following theorem which follows from the
global Torelli theorem.

\hfill

\theorem
Let $(M,I)$ be a hyperk\"ahler manifold, and
\[ \Mon(M,I)\subset O(H^2(M,\Z))\]
its monodromy group. Let $G$ the image of $\Aut(M)$
in $O(H^2(M,\Z))$. Then $G$ is the set of all 
$\gamma \in \Mon^{\Hdg}(M,I)$ fixing the K\"ahler chamber.

\hfill

{\bf Proof:} Similar to the
proof of \ref{_K3_auto_Theorem_} (second part); see \cite{_Markman:survey_}. 
\endproof

\hfill

Recall also that the image of the mapping class group is a finite index 
subgroup
in $O(H^2(M,\Z))$, and, accordingly, $\Mon^{\Hdg}(M,I)$ is of finite
index in the group of isometries of the Picard lattice.

\hfill

Let $M$ be a hyperk\"ahler manifold, and
$s\in H_2(M,\Z)$ a homology class. The 
BBF form defines an injection $H_2(M,\Z)\stackrel j\hookrightarrow H^2(M,\Q)$,
hence $q(s,s)$ can be rational. However, the denominators
of $\im j $ are divisors of the discriminant $\delta$ of $q$, 
hence $\delta j(s)$ is always integer, and
$\delta^2q(s,s)\in \Z$.

\hfill

\conjecture\label{_minim_curve_Conjecture_}
Let $M$ be a hyperk\"ahler manifold.
Then there exists a constant $C>0$, depending only on the deformation 
type of 
$M$, such that for any
primitive MBM class $s$ in $H^2(M,\Z)$
one has $|q(s,s)|<C$.

\hfill

\remark\label{mincurves-ht}
This conjecture is slightly weaker than its following, more 
algebraic-geometric version, implicitely appearing in \cite{bayer-hass-tschi}: let $M$
be a hyperk\"ahler manifold. Then there exists a constant $C>0$ such that
for any extremal rational curve (of minimal degree) $R$ on any deformation
$(M,I),\ I\in \Teich$, one has $|q(R,R)|<C$. 

Indeed, for any primitive
integral MBM class $s$, some integral multiple $Ns$ is represented by 
an extremal rational curve $R$ on such a deformation $(M,I)$ that
$s$ generates its Picard group, by \ref{_MBM_posi_theorem_}. 
So the boundedness of $|q(R,R)|$ means
the boundedness of $N$ plus the boundedness of $|q(s,s)|$. 

\hfill

\theorem\label{morkaw}
Let $M$ be a hyperk\"ahler manifold.
Then \ref{_minim_curve_Conjecture_} implies the
Morrison-Kawamata cone conjecture for all deformations of
$M$.

\hfill

{\bf Proof:}
Fix a complex structure $I$ on $M$, and let $S(I)$
be the set of MBM classes
which are of type $(1,1)$ on $(M,I)$. 
The faces of the K\"ahler-Weyl chambers are pieces of
$s^\bot$, where $s^\bot$ runs through $S(I)$.
If \ref{_minim_curve_Conjecture_} is true, the
monodromy group $\Mon^{\Hdg}(M, I)$ acts on $S(I)$
with finitely many orbits (\ref{_Kneser_orbits_Theorem_}).
Then the argument as in \ref{_cone_for_K3_Corollary_} 
proves that $\Mon^{\Hdg}(M, I)$ acts on the set of faces 
of K\"ahler-Weyl chambers with finitely many orbits.

Let ${\goth F}$ be the set of all pairs $(F,\nu)$, where 
$F$ is a face of a K\"ahler-Weyl chamber, and $\nu$ is 
orientation on a 
normal bundle $NF$. Then the monodromy acts on 
${\goth F}$ with finitely many orbits, as we have 
already indicated. Each face of a K\"ahler cone $\Kah$
gives an element of ${\goth F}$: we pick orientation
determined by the side of a face adjoint to the cone.

Denote by ${\goth F}_0\subset {\goth F}$ the set of 
all faces of the K\"ahler cone with their orientations.
There are finitely many orbits of $\Mon^{\Hdg}(M,I)$ acting on
${\goth F}_0$. However, each $\gamma\in \Mon^{\Hdg}(M,I)$
which maps an element $f\in {\goth F}_0$ to an 
element $\gamma(f) \in {\goth F}_0$ maps the K\"ahler chamber
to itself. Indeed, there are two K\"ahler-Weyl chambers
adjoint to each face, and $\gamma(\Kah)$ is one of the two
chambers adjoint to $\gamma(f)$. But since the
orientation is preserved, $\gamma(\Kah)=\Kah$ and not
the other one. Therefore $\gamma$ is induced by an automorphism
of $(M,I)$. \endproof

\hfill

For hyperk\"ahler manifolds which are deformation equivalent 
to the Hilbert scheme of length $n$ subschemes on a K3 surface, 
\ref{_minim_curve_Conjecture_} is easily deduced from \cite{bayer-hass-tschi}, 
Proposition 2. Indeed it is shown there that for $M$ as above
and projective, any extremal ray of the Mori cone contains
an effective curve class $R$ with $q(R,R)\geq -\frac{n+3}{2}$.
But if we want to bound the ``length'' (that is, the square of a
minimal representative curve) of an MBM class on $(M,I)$, we can
look at the extremal rays on projective deformations, because
MBM classes are deformation equivalent and the monodromy acts 
by isometries. We thus obtain the following 

\hfill

\corollary\label{morkaw-hilbK3}
The Morrison-Kawamata cone conjecture holds for deformations of
the Hilbert scheme of length $n$ subschemes on a K3 surface. \endproof

\hfill

{\bf Acknowledgements:}
We are
grateful to Eyal Markman for many interesting discussions
and an inspiration. Many thanks to Brendan Hassett and to Justin Sawon who disabused
us of some wrong assumptions. We are grateful to Claire Voisin
for patiently answering many questions, and to 
Fr\'ed\'eric Campana for his expert advice on cycle spaces of non-algebraic manifolds.
Many thanks to Wai Kiu Chan and William Jagy for answering our questions
on \url{http://mathoverflow.net} and for the reference to Kneser's book.

{

\noindent {\sc Ekaterina Amerik\\
{\sc Laboratory of Algebraic Geometry,\\
National Research University HSE,\\
Department of Mathematics, 7 Vavilova Str. Moscow, Russia,}\\
\tt  Ekaterina.Amerik@math.u-psud.fr}, also: \\
{\sc Universit\'e Paris-11,\\
Laboratoire de Math\'ematiques,\\
Campus d'Orsay, B\^atiment 425, 91405 Orsay, France}

\hfill

\noindent {\sc Misha Verbitsky\\
{\sc Laboratory of Algebraic Geometry,\\
National Research University HSE,\\
Department of Mathematics, 7 Vavilova Str. Moscow, Russia,}\\
\tt  verbit@mccme.ru}, also: \\
{\sc Kavli IPMU (WPI), the University of Tokyo}

}

\end{document}